\documentclass[11pt,twoside]{article}

\usepackage{amsfonts}
\usepackage[mathscr]{euscript}

\newtheorem{theo}{Theorem}[section]
\newtheorem{prop}{Proposition}[section]
\newtheorem{lem}{Lemma}[section]

\newcommand{\be}{\begin{equation}}
\newcommand{\ee}{\end{equation}}
\newcommand\bes{\begin{eqnarray}} \newcommand\ees{\end{eqnarray}}
\newcommand{\bess}{\begin{eqnarray*}}
\newcommand{\eess}{\end{eqnarray*}}
\newcommand\bedd{\bes\left\{\begin{array}{ll}}
\newcommand\eedd{\end{array}\right.\ees}

\newcommand\ep{\varepsilon}
\newcommand\kk{\left}
\newcommand\rr{\right}

\newcommand\dd{\displaystyle}

\begin{document}\thispagestyle{empty}

\begin{center}{\Large\bf A free boundary problem for a predator-prey model\\[2mm]
 with double free boundaries}\footnote{This work was supported by NSFC Grants 11071049 and 11371113}\\[5mm]
 {\Large  Mingxin Wang$^{a}$,\ \ Jingfu Zhao$^{a,b}$}\\[2mm]
{\small $a$. Natural Science Research Center, Harbin Institute of Technology, Harbin 150080, PR China}\\[1mm]
{\small $b$. Department of Mathematics, Shaanxi University of Technology, Hanzhong
723000, PR China. }
\end{center}

\begin{quote}\noindent{\bf Abstract.}


 \end{quote}

\def\theequation{\arabic{section}.\arabic{equation}}

\section{Introduction}
\setcounter{equation}{0}
{\setlength\arraycolsep{2pt}

The expanding (migration) of a new or invasive species is one of the most
important topics in mathematical ecology. A lot of mathematicians have made
efforts to develop various invasion models and investigated them from a viewpoint
of mathematical ecology. To describe the invasion and spreading phenomenon, there have been many interesting studies on the existence of positive traveling wave solutions connecting two different equilibria. Also, the study of asymptotic spreading speed plays an important role in invasion ecology since it can be used to predict the mean spreading rate of species. On the other hand, Du and Lin \cite{DLin} proposed a new
mathematical model to understand the expanding of an invasive or new species.
Their model is described as a free boundary problem for a logistic diffusion
equation:
 \be
 \left\{\begin{array}{lll}
 u_t-du_{xx}=u(a-bu), &t>0,\ \ 0<x<h(t),\\[1mm]
 u_x(t,0)=0,\ \ u(t,h(t))=0,\ \ &t>0,\\[1mm]
 h'(t)=-\mu u_x(t,h(t)),&t>0,\\[1mm]
 h(0)=h_0,\ \ u(0,x)=u_0(x),&0<x<h(t),
 \end{array}\right.\label{1.1}
 \ee
where $x=h(t)$ is the moving boundary to be determined, $a,b,d,\mu$ and $h_0$ are
given positive constants, and $u_0$ is a given positive initial function. The dynamics of the free boundary is determined by Stefan-like condition $h'(t)=-\mu u_x(t,h(t))$. This condition means that the population pressure at the free boundary is a driving force of the free boundary. Du and Lin \cite{DLin} have established the existence and uniqueness
of global solutions and, furthermore, derived various interesting results about the
long time behavior of solution. One of very remarkable results
is a spreading-vanishing dichotomy of the species, i.e., the solution $(u,h)$ of (\ref{1.1}) satisfies one of the following properties:

$\bullet$ {\bf Spreading}:\, $h(t)\to\infty,\,  u(t,x)\to a/b$ as $t\to\infty$;

$\bullet$ {\bf Vanishing}: \, $h(t)\to h_\infty\leq (\pi/2)\sqrt{d/a}$, and $u(t,x)\to 0$ as $t\to\infty$.\\
When the spreading occurs, it is also proved that the spreading speed approaches
to a positive constant $k_0 $, i.e., $h(t)=(k_0+\circ(1))t$ as $t\to\infty$. See
also the paper of Du and Guo \cite{DG,DG1}, where a free boundary problem similar to
(\ref{1.1}) was studied in higher space dimension and the same spreading-vanishing
dichotomy has been established. In \cite{KY}, (\ref{1.1}) was discussed with $u_x(t,0)=0$
replaced by $u(t,0)=0$.

A variety of reaction-diffusion systems are used to describe some phenomena arising in
population ecology. A typical model is the following classical Lotka-Volterra type predator-prey system in a one-dimensional habitat (under the suitable rescaling)
\bes
 \left\{\begin{array}{lll}
 u_t-u_{xx}=u(1-u+av), &t>0,\ \ x\in\mathbb{R},\\[1mm]
  v_t-Dv_{xx}=v(b-v-cu),\ \ &t>0, \ \ x\in\mathbb{R},
 \end{array}\right.\label{a.1}
 \ees
where $u(t,x), v(t,x)$ denote, respectively, the population densities of predaor and prey at the position $x$ and time $t$.

Understanding of spatial and temporal behaviors of interacting species in ecological systems is a central issue in population ecology. One aspect of great interest for a model with multispecies interactions is whether the species can spread successfully. Motivated by the work of Du and Lin \cite{DLin}, in the present paper we shall study a free boundary problem associated with (\ref{a.1}) to realize the expanding mechnism of the species. In the real world, the following two kind of phenomenons often occur:

(i)\, At the initial state, one kind of prey species (for example, {\it pest species}) occupied the whole space or a large region. In order to control such prey species we put one kind of predator species ({\it natural enemies}) in some bounded region or a small region (initial habitat).

(ii)\, There is some kind of species (prey) in the whole space or a large region, and at some time (initial time) another type species (the new or invasive species, predator) enters some bounded region or a small region (initial habitat).

In general, the predator has a tendency to emigrate from the boundaries to obtain their new habitat, i.e., it will move outward along the unknown curves (free boundaries) as time increases. It is assumed that the movement speeds of free boundaries are proportional to the gradient of predator. We want to realize the dynamics/variations of predator, prey and free boundaries. According to the above arguments, the model we are concerned here is the following free boundary problem
 \bes
 \left\{\begin{array}{lll}
 u_t-u_{xx}=u(1-u+av), &t>0,\ \ g(t)<x<h(t),\\[1mm]
 u(t,x)\equiv 0, &t\ge 0,\ \ x\not\in(g(t), h(t)),\\[1mm]
 v_t-Dv_{xx}=v(b-v-cu),\ \ &t>0, \ \ x\in\mathbb{R},\\[1mm]
 u=0,\ \ g'(t)=-\mu u_x, \ &t\ge 0, \ \ x=g(t),\\[1mm]
 u=0,\ \ h'(t)=-\mu u_x, \ &t\ge 0, \ \ x=h(t),\\[1mm]
 g(0)=-h_0,\ \ h(0)=h_0,&\\[2mm]
 u(0,x)=u_0(x), \ x\in[-h_0,h_0]; \ &
 v(0,x)=v_0(x), \ x\in\mathbb{R},
 \end{array}\right.\label{1.2}
 \ees
where $\mathbb{R}=(-\infty,\infty)$, $x=g(t)$ and $x=h(t)$ represent the left and
right moving boundaries, respectively, which are to be determined, $a, b, c, D, h_0$ and $\mu$ are given positive constants.  The initial functions $u_0(x),v_0(x)$ satisfy
 \bess
 u_0\in C^2([-h_0,h_0]),\ u_0(\pm h_0)=0,\ u_0>0\ \ {\rm in} \
(-h_0,h_0);\ \
 v_0\in C_b(\mathbb{R}),\ v_0>0\ \ {\rm in} \ \mathbb{R},
 \eess
here $C_b(\mathbb{R})$ is the space of continuous and bounded functions in $\mathbb{R}$. The ecological background of the free boundary conditions $g'(t)=-\mu u_x(t,g(t))$ and $h'(t)=-\mu u_x(t,h(t))$ can also refer to \cite{BDK}.

We will show that (\ref{1.2}) has a unique solution $(u(t,x),v(t,x),g(t),h(t))$ defined
for all $t>0$, with $u(t,x)\geq0$, $v(t,x)>0$, $g'(t)<0$ and $h'(t)>0$.
Moreover, a spreading-vanishing dichotomy holds for (\ref{1.2}),
namely, as time $t\to\infty$, either
 \begin{quote}
(i) \ the predator $u(t,x)$ successfully establishes itself in the new environment
(henceforth called spreading) in the sense that $g(t)\to-\infty$
and $h(t)\to \infty$. Moreover, both $u(t,x)$ and $v(t,x)$ go to positive constants for
the weakly hunting case $b>c$ and $ac<1$, while
$u(t,x)\to 1$, $v(t,x)\to 0$ for the strongly hunting case: $b\leq c$;\end{quote}
or
 \begin{quote}
(ii) \ the predator $u(t,x)$ fails to
establish and vanishes eventually (called vanishing), i.e., $h(t)-g(t)\to
h_\infty-g_\infty\leq \pi\sqrt{1/(1+ab)}$,
$\|u(t,x)\|_{C[g(t),h(t)]}\to 0$ and $v(t,x)\to b$.
 \end{quote}
The criteria for spreading and vanishing are the following: If the initial occupying area $[-h_0,h_0]$ is beyond a critical size, namely $2h_0\geq
\pi\sqrt{1/(1+ab)}$, then regardless of the initial population size $(u_0, v_0)$, spreading always happens. On the
other hand, if $2h_0<\pi\sqrt{1/(1+ab)}$, then whether spreading or
vanishing occurs is determined by the initial population size $(u_0,v_0)$ and the
coefficient $\mu$ in the Stefan condition.

In the absence of $v$, the problem (\ref{1.2}) is reduced to a one phase Stefan problem for the logitic model which has been systematically studied by many authors, see, for example \cite{BDK}, \cite{DG}--\cite{DLin}, \cite{DLou, KY, PZ} (including the higher dimension and heterogeneous environment case) and the references cited therein. The one phase Stefan free boundary condition in (\ref{1.2}) also arises in many other applications, for instance, in the modeling of wound healing \cite{ChenA}. As far as population models are concerned, \cite{Lin} used such a condition for a predator-prey system over a bounded interval, showing that the free boundary reaches the fixed boundary in finite time, and hence, the long-time dynamical behavior of the system is the same as the well-studied fixed boundary problem; and in \cite{MYS1}, a two-phase Stefan condition was used for a competition system over a bounded interval, where the free boundary separates the two competitors from each other in the interval. There is a vast literature on the Stefan problems, and some important theoretical advances can be found in \cite{CS, CR} and the references therein.

The other related works concerning free boundary problems for biological models, please refer to, for instance \cite{DMZ, GW, HIMN, HMS} and references cited therein.

The organization of this paper is as follows. In Section $2$, we first use a contraction
mapping argument to prove the local existence and uniqueness of solution to (\ref{1.2}), and then show that it exists for all time $t\in(0,\infty)$. In order to estimate $(u(t,x),v(t,x))$ and $(g(t),h(t))$, in Section 3 we give some comparison principles. Section 4 is devoted to the long time behavior of $(u(t,x),v(t,x))$. Theorem \ref{th4.2} plays key roles in the following two aspects: (i) affirming the predator species disappears eventually; (ii) determining the criteria for spreading and vanishing (see the following Section $5$). Moreover, its proof is very different from the single equation case (refer to the proofs of \cite[Lemma 3.1]{DLin} and \cite[Theorem 2.10]{KY}). In  Section 5 we shall provide the criteria for spreading and vanishing. The last section is a brief discussion.

Before ending this section, we should emphasize here that if $-h_0$ is replaced by another number $g_0$ with $g_0<h_0$, and/or the free boundary conditions $g'(t)=-\mu u_x(t,g(t))$ and $h'(t)=-\mu u_x(t,h(t))$ are replaced by $g'(t)=-\mu_1 u_x(t,g(t))$ and $h'(t)=-\mu_2 u_x(t,h(t))$, respectively, and $\mu_1,\,\mu_2$ are positive constants, then all results of the present paper are still true.

\section{Existence and uniqueness}
\setcounter{equation}{0}

In this section, we first prove the following local existence and uniqueness result by contraction mapping theorem, and then use suitable estimate to illustrate that the solution is defined for all $t>0$.

\begin{theo}\label{th2.1} \ For any given $\alpha\in(0,1)$, there is a $T>0$ such that problem $(\ref{1.2})$ admits a unique solution
$(u,v,g,h)\in C^{\frac{1+\alpha}2,1+\alpha}(\overline{D}_T)\times
\mathfrak{C}_T\times [C^{1+\frac \alpha 2}([0,T])]^2$. And
 \bes
\|u\|_{C^{\frac{1+\alpha}2,1+\alpha}(\overline{D}_T)}+\|g\|_{C^{1+\frac \alpha
2}([0,T])}+\|h\|_{C^{1+\frac \alpha 2}([0,T])}\leq C,\label{2.1}
  \ees
where
 \bess
 D_T=\kk\{0<t\le T,\, g(t)<x<h(t)\rr\},\ \ \mathfrak{C}_T=C_b([0,T]\times\mathbb{R})\cap C^{1+\frac{\alpha}2,2+\alpha}_{{\rm
 loc}}((0,T]\times\mathbb{R}),\eess
$C$ and
$T$ only depend on $h_0$, $\alpha$, $\|u_0\|_{W_p^2([-h_0,h_0])}$ with $p\geq(n+2)/(1-\alpha)$ and
$\|v_0\|_{C_b(\mathbb{R})}$.
 \end{theo}

{\bf Proof.}\, As in \cite{ChenA}, we first straighten the free boundaries.
Let $\zeta(y)$ be a function in $C^3(\mathbb{R})$ satisfying
 $$\zeta(y)=1\ \ \mbox{if}\ |y-h_0|<{h_0}/{4},\ \ \zeta(y)=0\ \ \mbox{if}\
|y-h_0|>{h_0}/{2},\  \ |\zeta'(y)|<{6}/{h_0}, \ \ \forall \ y\in\mathbb{R},$$
and set $\xi(y)=-\zeta(-y)$. Consider the transformation
 $$(t,x)\to(t,y),\ {\rm where}\ \ x=y+\zeta(y)(h(t)-h_0)
  +\xi(y)(g(t)+h_0),\ y\in\mathbb{R}.$$
Notice that as long as $|h(t)-h_0|+|g(t)+h_0|\leq{h_0}/{16}$, 
the above transformation is a diffeomorphism from $\mathbb{R}$ onto $\mathbb{R}$.
Moreover, it changes the free boundaries
$x=g(t)$, $x=h(t)$ to the lines $y=-h_0$ and $y=h_0$ respectively. Now, direct
calculations yield
   \bess
  &\dd\frac{\partial y}{\partial
x}=\dd\frac{1}{1+\zeta'(y)(h(t)-h_0)+\xi'(y)(g(t)+h_0)}
  \equiv\sqrt{\rho(g(t),h(t),y)},&\\[2mm]
  &\dd\frac{\partial^2 y}{\partial x^2}
  =-\dd\frac{\zeta''(y)(h(t)-h_0)+\xi''(y)(g(t)+h_0)}
   {[1+\zeta'(y)(h(t)-h_0)+\xi'(y)(g(t)+h_0)]^3}\equiv \varrho(g(t),h(t),y),&\\[2mm]
  &\dd\frac{\partial y}{\partial t}=-\dd\frac{\zeta(y)h'(t)
  +\xi(y)g'(t)}{1+\zeta'(y)(h(t)-h_0)+\xi'(y)(g(t)+h_0)}
   \equiv \varsigma(g(t),g'(t),h(t),h'(t),y).&\eess

If we set
 $$\begin{array}{l}
u(t,x)=u(t,y+\zeta(y)(h(t)-h_0)+\xi(y)(g(t)+h_0))=w(t,y),\\[2mm]
 v(t,x)=v(t,y+\zeta(y)(h(t)-h_0)+\xi(y)(g(t)+h_0))=z(t,y),
 \end{array}$$
then
 \[\begin{array}{lll}
 u_t=w_t+\varsigma(g(t),g'(t),h(t),h'(t),y)w_y,
&v_t=z_t+\varsigma(g(t),g'(t),h(t),h'(t),y)z_y,\\[2mm]
 u_x=\sqrt{\rho(g(t),h(t),y)}w_y, &v_x=\sqrt{\rho(g(t),h(t),y)}z_y,\\[2mm]
 u_{xx}=\rho(g(t),h(t),y)w_{yy}+\varrho(g(t),h(t),y)w_y, \ \ &
v_{xx}=\rho(g(t),h(t),y)z_{yy}+\varrho(g(t),h(t),y)z_y
 \end{array}\]
and $(w,z)$ satisfies
 \bes
 \left\{\begin{array}{lll}
 w_t-\rho w_{yy}-(\varrho -\varsigma)w_y=w(1-w+az),\ &t>0, \
|y|<h_0,\\[1mm]
 w(t,y)\equiv 0, &t\ge 0, \ |y|\geq h_0,\\[1mm]
 z_t-D\rho z_{yy}-(D\varrho-\varsigma)z_y=z(b-z-cw), \ &t>0,\ \
y\in\mathbb{R},\\[1mm]
 w(t,-h_0)=w(t,h_0)=0,&t\ge 0,\\[1mm]
 w(0,y)=u_0(y), \ y\in[-h_0,h_0]; \ 
 z(0,y)=v_0(y), \ &y\in\mathbb{R},
 \end{array}\right.\label{2.2}
 \ees
where $\rho=\rho(g(t),h(t),y)$, $\varrho=\varrho(g(t),h(t),y)$,
$\varsigma=\varsigma(g(t),g'(t),h(t),h'(t),y)$.

\vskip 2pt
Let $g^*=-\mu u_0'(-h_0)$ and $h^*=-\mu u_0'(h_0)$. For $0<T<\frac{h_0}{16}\min\big\{\frac 1{1+g^*},\,\frac 1{1+h^*}\big\}$, 
we define $I_T=[0,T]\times[-h_0,h_0]$, and
 $$ \begin{array}{l}
 \mathcal{D}^1_T=\big\{w\in C^{\frac \alpha 2,\alpha}(I_T):\
 w(t,\pm h_0)=0,\, w(0,y)=u_0(y),\,  
  \|w-u_0\|_{C^{\frac \alpha 2,\alpha}(I_T)}\leq 1\big\},\\[2mm]
\mathcal{D}^2_T=\{g\in C^1([0,T]):\, g(0)=-h_0,\ g'(0)=g^*, \
\|g'-g^*\|_{C([0,T])}\leq 1\},\\[2mm]
 \mathcal{D}^3_T=\{h\in C^1([0,T]):\, h(0)=h_0,\ h'(0)=h^*, \
\|h'-h^*\|_{C([0,T])}\leq 1\}.
 \end{array}
 $$
It is easily seen that $\mathcal{D}_T= \mathcal{D}^1_T\times \mathcal{D}^2_T\times
\mathcal{D}^3_T$
 is a closed convex set in $C^{\frac \alpha 2,\alpha}(I_T)\times C^1([0,T])\times
C^1([0,T])$.

Next, we shall prove the existence and uniqueness result by using the contraction
mapping theorem. First, we observe that due to our choice of $T$, for any given
$(w,g,h)\in \mathcal{D}_T$, there holds:
 $$ |g(t)+h_0|\leq T\|g'\|_{C([0,T])}\leq {h_0}/{16},\ |h(t)-h_0|\leq
 T\|h'\|_{C([0,T])}\leq {h_0}/{16}.$$
Therefore the transformation $(t,y)\to(t,x)$ introduced at the beginning of the proof
is well defined. For any $(w,g,h)\in \mathcal{D}_T$, let $\hat w(t,y)=w(t,y)$ when
$|y|\leq h_0$, and $\hat w(t,y)=0$ when $|y|>h_0$. Since $w\in C^{\frac \alpha
2,\alpha}(I_T)$, it is easy to see that $\hat w\in C^{\frac \alpha
2,\alpha}([0,T]\times\mathbb{R})$. The standard partial differential equation theory
\cite{Fri,LSU} guarantees that the problem
 \bess
 \left\{\begin{array}{lll}
 z_t-D\rho z_{yy}-(D\varrho-\varsigma) z_y=z(b-z-c\hat w),\ &t>0,\ \
y\in\mathbb{R},\\[1mm]
 z(0,y)=v_0(y),&y\in\mathbb{R}
 \end{array}\right.\eess
admits a unique solution $z\in C_b([0,T]\times\mathbb{R})\cap C^{1+\frac\alpha
2,2+\alpha}_{\rm loc}((0,T]\times\mathbb{R})$. Also, the following initial boundary
value problem
 \[ \left\{\begin{array}{lll}
 \tilde w_t-\rho \tilde w_{yy}-(\varrho -\varsigma)\tilde w_y=w(1-w+az),
\ &t>0, \ -h_0<y<h_0,\\[1mm]
 \tilde w(t,-h_0)=\tilde w(t,h_0)=0,&t>0,\\[1mm]
 \tilde w(0,y)=u_0(y),&-h_0<y<h_0
 \end{array}\right.\]
admits a unique solution $\tilde w\in C^{\frac{1+\alpha}2,1+\alpha}(I_T)$.
Moreover, using $L^p$ estimate for parabolic equations with $p\geq(n+2)/(1-\alpha)$ and Sobolev's inequalities, one gets
 \bes
 \|\tilde w\|_{C^{\frac{1+\alpha}2,1+\alpha}(I_T)}\leq C_1.\label{2.3}
 \ees
where $C_1$ is a constant dependent on $h_0$, $\alpha$, $\|u_0\|_{W_p^2([-h_0,h_0])}$ and $\|v_0\|_{C_b(\mathbb{R})}$. Define
 $$ \tilde g(t)=-h_0-\int_0^t\mu\tilde w_y(\tau,-h_0){\rm d}\tau,\ \ \tilde
 h(t)=h_0-\int_0^t\mu\tilde w_y(\tau,h_0){\rm d}\tau. $$
Then $\tilde g'(t)=-\mu\tilde w_y(t,-h_0), \tilde h'(t)=-\mu\tilde w_y(t,h_0)$. Subsequently, $\tilde g'(t), \tilde h'(t)\in C^{\frac \alpha 2}([0,T])$,
and
 \bes
\|\tilde g'(t)\|_{C^{\frac \alpha 2}([0,T])},\ \|\tilde h'(t)\|_{C^{\frac \alpha
2}([0,T])}\leq \mu C_1:=C_2.\label{2.4}
 \ees

We now define $\mathcal{F}: \mathcal{D}_T\rightarrow C^{\frac \alpha
2,\alpha}(I_T)\times C^1([0,T])\times C^1([0,T])$ by
 \[\mathcal{F}(w,g,h)=(\tilde w,\tilde g,\tilde h).\]
Clearly $(w,g,h)\in \mathcal{D}_T$ is a fixed point of $\mathcal{F}$ if and only if
$(w,z,g,h)$ solves (\ref{2.2}). By (\ref{2.3}) and (\ref{2.4}), one has
 \bess
\|\tilde g'-g^*\|_{C([0,T])}\leq\|\tilde g'\|_{C^{\frac\alpha
2}([0,T])}T^{\frac\alpha 2}\leq \mu C_1T^{\frac \alpha 2},\ \ 
\|\tilde h'-h^*\|_{C([0,T])}\leq\|\tilde h'\|_{C^{\frac\alpha
2}([0,T])}T^{\frac\alpha 2}\leq \mu C_1T^{\frac \alpha 2},
 \eess
and
  \bess
\|\tilde w-u_0\|_{C^{\frac \alpha 2,\alpha}(I_T)}
&=&\|\tilde w-u_0\|_{C(I_T)}+[\tilde w-u_0]_{C^{\frac \alpha
2,\alpha}(I_T)}\\[.5mm]
&\leq &\|\tilde w\|_{C^{\frac{1+\alpha}2,0}(I_T)}T^{\frac{1+\alpha}2}+\|\tilde
w\|_{C^{\frac{1+\alpha}2,0}(I_T)}T^{\frac 1 2}+
(2h_0)^{1-\alpha}\|\tilde w_y\|_{C^{\frac \alpha 2,0}(I_T)}T^{\frac \alpha
2}\\[.5mm]
&\leq& C_1\kk(T^{\frac{1+\alpha}2}+T^{\frac 1 2}+(2h_0)^{1-\alpha}T^{\frac\alpha
2}\rr).
  \eess
Therefore, if we take $T\leq\min\big\{1,(\mu C_1)^{-2/\alpha},\
[(2+(2h_0)^{1-\alpha})C_1]^{-2/\alpha}\big\} $, then $\mathcal{F}$ maps
$\mathcal{D}_T$ into itself.

Next we attest that $\mathcal{F}$ is a contraction mapping on $\mathcal{D}_T$ for $T>0$
sufficiently small. Indeed, let $(w_i,g_i,h_i)\in \mathcal{D}_T\ (i=1,2)$ and denote $(\tilde w_i,\tilde
g_i,\tilde h_i)=\mathcal{F}(w_i,g_i,h_i)$.
Then it follows from (\ref{2.3}) and (\ref{2.4}) that
  $$ \|\tilde w_i\|_{C^{\frac{1+\alpha}{2},1+\alpha}(I_T)}\leq C_1,\
 \|\tilde g_i'(t)\|_{C^{\frac \alpha 2}([0,T])}\leq C_2,\ \|\tilde h_i'(t)\|_{C^{\frac
\alpha 2}([0,T])}\leq C_2.$$
Setting $\gamma=\tilde w_1-\tilde w_2,\, \zeta=z_1-z_2$, we find that $\gamma$ and $\zeta$
satisfy
  \bess\left\{\begin{array}{lll}
\gamma_t-\rho_1\gamma_{yy}-(\varrho_1-\varsigma_1)\gamma_y=(\rho_1-\rho_2)\tilde w_{2yy}+(\varrho_1-\varrho_2)\tilde
w_{2y}+(\varsigma_2-\varsigma_1)\tilde w_{2y}\\[1mm]
 \qquad\qquad +[1-(w_1+w_2)-az_1](w_1-w_2)+aw_2\zeta,\ \ t>0,\ -h_0<y<h_0,\\[1mm]
\gamma(t,\pm h_0)=0,\ t\ge 0; \ \ \gamma(0,y)=0,\ \ -h_0\leq y\leq h_0
  \end{array}\right.\eess
and
 \bess\left\{\begin{array}{lll}
\zeta_t-D\rho_1\zeta_{yy}-(D\varrho_1-\varsigma_1)\zeta_y-[b-c\hat w_1-(z_1+z_2)]\zeta\\[1mm]
\ \ \
=D(\rho_1-\rho_2)z_{2yy}+D(\varrho_1-\varrho_2)z_{2y}+(\varsigma_2-\varsigma_1)z_{2y}-c(\hat
w_1-\hat w_2)z_2,\ t>0,\ y\in\mathbb{R},\\[1mm]
\zeta(0,y)=0,\ \ y\in\mathbb{R},
 \end{array}\right.\eess
respectively, where $\rho_i=\rho(g_i(t),h_i(t),y)$, $\varrho_i=\varrho(g_i(t),h_i(t),y)$ and $\varsigma_i=\varsigma(g_i(t),g_i'(t),h_i(t),h_i'(t),y)$, 
$i=1,2$. In view of the standard theory for parabolic partial differential equations and Sobolev's imbedding theorem \cite{LSU}, we obtain that
 \bes\|z_1-z_2\|_{C(I_T)}&\leq&
 C_3\dd\kk(\|w_1-w_2\|_{C(I_T)}+\|(g_1-g_2,\, h_1-h_2)\|_{C^1([0,T])}\rr),
 \nonumber\\[1mm]
\|\tilde w_1-\tilde w_2\|_{C^{\frac{1+\alpha}2,1+\alpha}(I_T)}
&\leq &C_3\dd\kk(\|(w_1-w_2,\, z_1-z_2)\|_{C(I_T)}+\|(g_1-g_2,\,h_1-h_2)\|_{C^1([0,T])}\rr)\nonumber\\[1mm]
 &\leq &C_4\dd\kk(\|w_1-w_2\|_{C(I_T)}+\|(g_1-g_2,\,h_1-h_2)\|_{C^1([0,T])}\rr),\qquad
 \label{2.5}\ees
where $C_3$ and $C_4$ depend on $C_1, C_2$ and the functions $\rho, \varrho, \varsigma$.
Taking the difference of equations for
$\tilde g_1,\tilde h_1,\tilde g_2,\tilde h_2$ results in
 \be
\|\tilde g_1'-\tilde g_2'\|_{C^{\frac \alpha 2}([0,T])},\ \|\tilde h_1'-\tilde
h_2'\|_{C^{\frac \alpha 2}([0,T])}
\leq \mu\|\tilde w_{1y}-\tilde w_{2y}\|_{C^{\frac \alpha 2,0}(I_T)}.\label{2.6}
 \ee
We may assume that $T\leq 1$. Combining (\ref{2.5}) and (\ref{2.6}), and applying the mean value theorem, it yields
 \bess
&&\|\tilde w_1-\tilde w_2\|_{C^{\frac{1+\alpha}2,1+\alpha}(I_T)}+\|(\tilde
g_1'-\tilde g_2',\,\tilde h_1'-\tilde h_2')\|_{C^{\frac \alpha 2}([0,T])}\\[1mm]
&\leq & C_5
\kk(\|w_1-w_2\|_{C(I_T)}+\|(g_1'-g_2',\,h_1'-h_2')\|_{C([0,T])}\rr),
  \eess
where $C_5$ depends on $C_4$ and $\mu$. On the other hand, by direct calculations,
 \bess
 \|\tilde w_1-\tilde w_2\|_{C^{\frac \alpha 2,\alpha}(I_T)}&\leq&\|\tilde
w_1-\tilde w_2\|_{C^{\frac{1+\alpha}2,0}(I_T)}T^{\frac{1+\alpha}2}
 +\|\tilde w_1-\tilde w_2\|_{C^{\frac{1+\alpha}2,0}(I_T)}T^{\frac 1 2}\\[1mm]
&&+(2h_0)^{1-\alpha}\|\tilde w_{1y}-\tilde w_{2y}\|_{C^{\frac \alpha
2,0}(I_T)}T^{\frac \alpha 2}\\[1mm]
&\leq &\kk(2+(2h_0)^{1-\alpha}\rr)T^{\frac \alpha 2}\|\tilde w_1-\tilde
w_2\|_{C^{\frac{1+\alpha}2,1+\alpha}(I_T)}.
  \eess
Let $\ep_1=h_0/16$ and $\ep_2=2+(2h_0)^{1-\alpha}$. Then for
 \[T:=\min\kk\{1,\ \, \frac{\ep_1}{1+g^*},\ \, \frac{\ep_1}{1+h^*},\, \ \frac 1{(\mu
 C_1)^{\frac 2\alpha}},\ \,
 \frac 1{(\ep_2C_1)^{2/\alpha}},\ \, \frac 1{(2\ep_2C_5)^{ 2/\alpha}}\rr\},\]
it follows that
 \bess
&&\dd\|\tilde w_1-\tilde w_2\|_{C^{\frac \alpha 2,\alpha}(I_T)}+\|(\tilde g_1'-\tilde
g_2',\,\tilde h_1'-\tilde h_2')\|_{C([0,T])}\\[1mm]
&\leq &\dd\kk(2+(2h_0)^{1-\alpha}\rr)T^{\frac \alpha 2}\|\tilde w_1-\tilde
w_2\|_{C^{\frac{1+\alpha}2,1+\alpha}(I_T)}
+T^{\frac \alpha 2}\|(\tilde g_1'-\tilde g_2',\,\tilde h_1'-\tilde h_2')\|_{C^{\frac \alpha
2}([0,T])}\\[1mm]
 &\leq &\dd(2+(2h_0)^{1-\alpha})C_5T^{\frac \alpha
2}\kk(\|w_1-w_2\|_{C(I_T)}+\|(g_1'-g_2',\,h_1'-h_2')\|_{C([0,T])}\rr)\\[1mm]
&\leq &\dd\frac{1}{2}\kk(\|w_1-w_2\|_{C^{\frac \alpha
2,\alpha}(I_T)}+\|(g_1'-g_2',\,h_1'-h_2')\|_{C([0,T])}\rr).
 \eess

The above arguments ensure that the operator $\mathcal{F}$ is contractive on
$\mathcal{D}_T$. It now follows from the contraction mapping theorem that
$\mathcal{F}$ has a unique fixed point $(w,g,h)$ in $\mathcal{D}_T$. Moreover,
by the $L^p$ estimates, we have additional regularity for $(w,z,g,h)$ as a solution of
(\ref{2.2}), namely, $w\in C^{\frac{1+\alpha}2,1+\alpha}(I_T)$,
$z\in C_b([0,T]\times\mathbb{R})\cap C^{1+\frac\alpha 2,2+\alpha}_{\rm
loc}((0,T]\times\mathbb{R})$, and $g,\,h\in C^{1+\frac{\alpha}2}([0,T])$,
and (\ref{2.3}), (\ref{2.4}) hold. In other words, $(w,z,g,h)$ is the unique
local classical solution of the problem (\ref{2.2}). Hence, $(u,v,g,h)$ is
the unique classical solution of (\ref{1.2}). \ \ \ \ \fbox{}

To show that the local solution obtained in Theorem \ref{th2.1}
can be extended to all $t>0$, we need the following estimate.

\begin{lem}\label{lm2.1} \
The solution of the free boundary problem $(\ref{1.2})$ satisfies
  \bess
 &0<u(t,x)\leq M_1,\ \ 0< t\leq T,\ g(t)<x<h(t),&\\[1mm]
&0<v(t,x)\leq M_2,\  \ 0< t\leq T,\ x\in\mathbb{R},&\\[1mm]
&-M_3\leq g'(t)<0,\ 0<h'(t)\leq M_3,\ \ 0<t\leq T,&
  \eess
where $M_i$ is independent of $T$ for $i=1,2,3$.
\end{lem}

{\bf Proof.} Using the strong maximum principle, we are easy to see that $u>0$ in
$(0,T]\times(g(t),h(t))$ and $v>0$ in $(0,T]\times \mathbb{R}$ as long as the solution exists. Since $v(t,x)$ satisfies
  $$\left\{\begin{array}{ll}
v_t-Dv_{xx}=v(b-v-cu), \ &x\in\mathbb{R}, \ t>0,\\[1mm]
v(0,x)=v_0(x)> 0,&x\in\mathbb{R},
\end{array}\right.$$
it is obvious that $v\leq \max\kk\{\|v_0\|_\infty,\, b\rr\}:= M_1$. Similarly, as
$u$ satisfies
 $$ \left\{\begin{array}{lll}
 u_t-u_{xx}=u(1-u+av),\ \ &t>0,\ \ g(t)<x<h(t),\\[1mm]
 u(t,g(t))=u(t,h(t))=0,\ \ \ &t>0,\\[1mm]
  u(0,x)=u_0(x)> 0,&-h_0<x<h_0,
 \end{array}\right.$$
we also have $u\leq \max\kk\{\|u_0\|_\infty,\,(1+aM_1)\rr\}:= M_2$.

To prove $h'(t)>0$ for $0<t\leq T$, we use the transformation
  \[
 y=x/h(t),\  \ \ w(t,y)=u(t,x),\ \ \ z(t,y)=v(t,x)\]
to straighten the free boundary $x=h(t)$. A series of detailed calculation asserts
 \bess
\left\{\begin{array}{ll}
w_t-f(t)w_{yy}-\phi(t,y)w_y=w(1-w+az), \ &0<t\leq T,\ 0<y<1,\\[1mm]
w(t,0)>0, \ \ w(t,1)=0, \ \ &0<t\leq T,\\[1mm]
w(0,y)=u_0(h_0y),&0\leq y\leq 1,
\end{array}\right.
 \eess
where $f(t)=h^{-2}(t)$, $\phi(t,y)=yh'(t)/h(t)$. This is an initial-boundary value problem with fixed boundary. Since $w(t,y)>0$ for $t>0$ and $0\leq y<1$, by the Hopf boundary lemma, we have $w_y(t,1)<0$ for $t>0$. This combines with the relation $u_x=h^{-1}(t)w_y$ yields $u_x(t,h(t))<0$, and so $h'(t)>0$ for $t>0$. Similarly, $g'(t)<0$ for $t>0$.

Now we illustrate that $g'(t)\geq -M_3$ and $h'(t)\leq M_3$ for all $t\in(0,T)$ with
some $M_3$ independent of $T$. To this aim, let $M$ be a positive constant,
$ \Omega_M=\kk\{0<t<T,\ g(t)<x<g(t)+1/M\rr\}$, 
and construct an auxiliary function
  $$w(t,x)=M_2[2M(x-g(t))-M^2(x-g(t))^2].$$
We will choose $M$ so that $w\geq u$ in $\Omega_M$.

Direct calculations indicate that, for $(t,x)\in\Omega_M$,
 \bess
 &w_t=2M_2Mg'(t)(-1-M(g(t)-x))\geq 0,&\\[1mm]
 &-w_{xx}=2M_2M^2,\ \ \ u(1-u+av)\leq M_2(1+aM_1).&
 \eess
Therefore,
 $$w_t-w_{xx}\geq 2M_2M^2\geq M_2(1+aM_1)\geq u(1-u+av)\ \ {\rm in}\ \Omega_M$$
provided  $M^2\geq (1+aM_1)/2$. It is obvious that
  $$w(t,g(t)+ M^{-1})=M_2\geq u(t,g(t)+M^{-1}),\ \ w(t,g(t))=0=u(t,g(t)).$$
Because
  \bess
  &u_0(x)=\dd\int_{-h_0}^xu_0'(y)dy\leq (x+h_0)\|u_0'\|_{C[-h_0,h_0]}\ \ \
  {\rm on}\ \ [-h_0,-h_0+M^{-1}],&\\[1mm]
 &w(0,x)=M_2[2M(x+h_0)-M^2(x+h_0)^2]\geq M_2M(h_0+x)\ \ \
  {\rm on} \ \ [-h_0,-h_0+M^{-1}],&\eess
it is easily to see that if $MM_2\geq \|u_0'\|_{C[-h_0,h_0]}$, then
  $$u_0(x)\leq (x+h_0)\|u_0'\|_{C[-h_0,h_0]}\leq w(0,x) \ \ \ {\rm on } \ \
  [-h_0,-h_0+M^{-1}].$$
Let
 $$M=\max\kk\{\sqrt{\dd\frac{1+aM_1}{2}},
  \ \dd\frac{\|u_0'\|_{C[-h_0,h_0]}}{M_2}\rr\}.$$
We can apply the maximum principle to $w-u$ over $\Omega_M$ and deduce that $u(t,x)\leq
w(t,x)$ for $(t,x)\in\Omega_M$. It would then follow that $u_x(t,g(t))\leq w_x(t,g(t))=2M_2M$, and hence
 $$g'(t)=-\mu u_x(t,g(t))\geq -2\mu M_2M:=-M_3.$$

Similarly, we can proved $h'(t)\leq M_3$ for $0<t<T$. 
The proof is complete. \ \ \ \ \fbox{}

\begin{theo}\label{th2.2} \ The solution of problem $(\ref{1.2})$ exists and is unique for
all $t\in (0,\infty)$.
 \end{theo}

{\bf Proof.} Let $[0,T_{\rm max})$ be the maximal time interval in which the
solution exists. By Theorem \ref{th2.1}, $T_{\rm max}>0$. It remains to show that
$T_{\rm max}=\infty$. Arguing indirectly, it is assumed that $T_{\rm max}<\infty$. By
Lemma \ref{lm2.1}, there exist positive constants $M_1$,  $M_2$ and $M_3$, independent of
$T_{\rm max}$, such that
 \bess
  &0\leq u(t,x)\leq M_1 \ \ {\rm in} \ \ [0,T_{\rm max})\times[g(t),h(t)],\ \ 0\leq
v(t,x)\leq M_2\ \ {\rm in}\ \ [0,T_{\rm max})\times\mathbb{R},\\[1mm]
  &-M_3\leq g'(t)\leq 0,\ -M_3t\leq g(t)+h_0\leq 0,\ 0\leq h'(t)\leq M_3,\ 0\leq
 h(t)-h_0\leq M_3t\ \ {\rm in}\ \ [0,T_{\rm max}).&\eess
We now fix $\delta\in[0,T_{\rm max})$. From the proof
of Theorem \ref{th2.1}, it is easily seen that $v(t,\cdot)\in C_{{\rm loc}}^{1+\alpha/2}(\mathbb{R})$
and $\|v(t,\cdot)\|_{C(\mathbb{R})}\leq M_2$ for $t\in[\delta,T_{\rm max})$.
For $p\geq(n+2)/(1-\alpha)$, by the standard $L^p$ estimates and the Sobolev embedding theorem, we can find a constant
$C>0$, depending only on $\delta$ and $M_i$ $(i=1,2,3)$, such that
$\|u(t,\cdot)\|_{W^2_p([g(t),h(t)])}\leq C$ for $t\in[\delta,T_{\rm max})$. It then
follows from the proof of Theorem \ref{th2.1} that there exists a $\tau>0$,
depending only on $C$ and $M_i$ $(i=1,2,3)$, such that the solution of
(\ref{1.2}) with initial time $T_{\rm max}-\tau/2$ can be extended uniquely to the
time $T_{\rm max}-\tau/2+\tau$. But this contradicts the assumption. The proof is now
complete.\ \ \ \ \fbox{}

\section{Comparison principles}
\setcounter{equation}{0}

In this section we shall give some comparison principles which can be used to estimate the solution $(u(t,x),v(t,x))$ and the free boundaries $x=g(t)$ and $x=h(t)$.

\begin{lem} $($Comparison principle$)$\label{lm3.1} \
Let $T>0$, $\bar g,\bar h\in C^1([0,T])$ and $\bar g<\bar h$ in $[0,T]$. Let $\bar u\in C(\overline{O})\cap C^{1,2}(O)$ with $O=\{0<t\leq T,\, 
\bar g(t)<x<\bar h(t)\}$, and $\bar v\in C^{1,2}(G)$ with $G=(0,T]\times\mathbb{R}$. Assume that $(\bar u,\bar v,
\bar g,\bar h)$ satisfies
$$
 \left\{\begin{array}{ll}
  \bar u_t-\bar u_{xx}\geq\bar u(1-\bar u+a\bar v),\ \ &0<t\leq T,\ \ \bar
g(t)<x<\bar h(t),\\[1mm]
 \bar v_t-D\bar v_{xx}\geq\bar v(b-\bar v),&0<t\leq T,\ \ x\in\mathbb{R},\\[1mm]
 \bar u(t,\bar g(t))=0,\ \ \bar g'(t)\leq-\mu \bar u_x(t,g(t)),&0<t\leq T,\\[1mm]
 \bar u(t,\bar h(t))=0,\ \ \bar h'(t)\geq-\mu\bar u_x(t,h(t)),\ \ &0<t\leq T.
 \end{array}\right.$$
If $\bar g(0)\leq-h_0,\,\bar h(0)\geq h_0$, $\bar u(0,x)\geq 0$ on $[\bar g(0),\bar h(0)]$, $u_0(x)\leq\bar u(0,x)$ on $[-h_0,h_0]$ and $v_0(x)\leq\bar v(0,x)$ in $\mathbb{R}$,
then the solution $(u,v,g,h)$ of $(\ref{1.2})$ satisfies
  \bess
  g(t)\geq\bar g(t),\ h(t)\leq\bar h(t) \ \ {\rm on}\  [0,T]; \ \ 
 u(t,x)\leq\bar u(t,x)\ \  {\rm on}\  \overline{D}_T; \ \
  v(t,x)\leq\bar v(t,x)\ \  {\rm on}\  \overline{G},\eess
where $D_T$ is defined as in Theorem $\ref{th2.1}$.
\end{lem}

{\bf Proof.} For small $\varepsilon>0$, let
$(u_\varepsilon,v_\varepsilon,g_\varepsilon,h_\varepsilon)$ be the unique solution
of (\ref{1.2}) with $h_0,\mu, u_0(x)$ and $v_0(x)$ replaced by
$h_0^\varepsilon=(1-\varepsilon)h_0$, $\mu_\varepsilon=(1-\varepsilon)\mu$,
$u_0^\varepsilon(x)\in C^2([-h_0^\varepsilon,h_0^\varepsilon])$ and
$v_0^\varepsilon(x)\in C(\mathbb{R})$, respectively, where $u_0^\varepsilon(x)$ and
$v_0^\varepsilon(x)$ satisfy $u_0^\varepsilon(-h_0^\varepsilon)=u_0^\varepsilon(h_0^\varepsilon)=0$ and
 $$ 0< u_0^\varepsilon(x)\leq u_0(x)\ \ \, {\rm on}\ \ [-h_0^\varepsilon,h_0^\varepsilon],
 \ \ 0<v_0^\varepsilon(x)<v_0(x)\ \ \, {\rm in}\ \ \mathbb{R},$$
and as $\varepsilon\to 0$,
  $$u_0^\varepsilon\kk({h_0\over h_0^\varepsilon}x\rr)\to u_0(x) \ \ \,
  {\rm in}\ \ C^2([-h_0,h_0]),\ \ \
  v_0^\varepsilon\kk({h_0\over h_0^\varepsilon}x\rr)\to v_0(x)\ \ \,
 {\rm in}\ \ C(\mathbb{R}).$$
Apply the comparison principle to $\bar v$ and $v_\varepsilon$ we have
$v_\varepsilon<\bar v$ on $[0,T]\times\mathbb{R}$.

We claim that $g_\varepsilon(t)>\bar g(t)$ and $h_\varepsilon(t)<\bar h(t)$ for all
$t\in(0,T]$. Clearly, this is true for small $t>0$. If our claim does not hold,
then we can find a first $\tau\leq T$ such that $g_\varepsilon(t)>\bar g(t)$ and
$h_\varepsilon(t)<\bar h(t)$ for all $t\in(0,\tau)$
and at least one of $g_\varepsilon(\tau)=\bar g(\tau)$ and $h_\varepsilon(\tau)=\bar
h(\tau)$ is true. Without loss of generality,
it is assumed that $h_\varepsilon(\tau)=\bar h(\tau)$. Then
  \be h'_\varepsilon(\tau)\geq\bar h'(\tau).\label{3.1}\ee
For any $0<\sigma\leq T$, let
 \bes
 D_{\sigma}^\ep:=\{(t,x)\in \mathbb{R}^2: 0<t\leq\sigma,\
 g_\varepsilon(t)<x<h_\varepsilon(t)\}.\label{3.2}\ees
The strong maximum principle yields $u_\varepsilon(t,x)<\bar u(t,x)$ in $D^\ep_{\tau}$.
Obviously, $\bar u(t,h_\varepsilon(\tau))=u_\varepsilon(t,h_\varepsilon(\tau))=0$ since $h_\varepsilon(\tau)=\bar h(\tau)$. Therefore,
$\bar u_x(t,h_\varepsilon(\tau))\leq u_{\varepsilon x}(t,h_\varepsilon(\tau))$. Note
that
$u_{\varepsilon x}(\tau,h_\varepsilon(\tau))<0$ and $\mu_\varepsilon<\mu$, it follows
that $h'_\varepsilon(\tau)<\bar h'(\tau)$.
This contradicts to (\ref{3.1}). So, $g_\varepsilon(t)>\bar g(t)$ and
$h_\varepsilon(t)<\bar h(t)$ for all $t\in(0,T]$.
We may now apply the usual comparison principle over $D^\ep_T$ to conclude that
$u_\varepsilon<\bar u$
in $D^\ep_T$.

Since the unique solution of (\ref{1.2}) depends continuously on the parameters in
(\ref{1.2}), as $\varepsilon\to 0$,
$(u_\varepsilon,v_\varepsilon,g_\varepsilon,h_\varepsilon)$ converges to $(u,v,g,h)$,
the unique solution of (\ref{1.2}). The desired
result then follows by letting $\varepsilon\to 0$ in the inequalities
$u_\varepsilon<\bar u$, $v_\varepsilon<\bar v$,
$g_\varepsilon<\bar g$ and $h_\varepsilon<\bar h$. \ \ \ \ \fbox{}

\vskip 4pt The pair $(\bar u,\bar v,\bar g,\bar h)$ in Lemma \ref{lm3.1} is usually called an upper solution of (\ref{1.2}). In the same way as the proof of Lemma \ref{lm3.1}, we can prove the following two lemmas: 

\begin{lem} $($Comparison principle$)$\label{lm3.2}\, Let $T>0$, $\tilde h\in C^1([0,T])$ with $\tilde h>0$ on
$[0,T]$, $\tilde u\in C(\overline{\cal D})\cap C^{1,2}({\cal D})$ with
${\cal D}=\{0<t\leq T,\ -h_0<x<\tilde h(t)\}$.
Assume that  $(\tilde u, \tilde h)$ satisfies $\tilde h(0)\leq h_0$, $u_0(x)\geq\tilde u(0,x)$ on $[-h_0,\tilde h(0)]$, and 
 \bess
 \left\{\begin{array}{ll}
  \tilde u_t-\tilde u_{xx}\leq\tilde u(1-\tilde u),\ \ &0<t\leq T,\, \ -h_0<x<\tilde
h(t),\\[1mm]
  \tilde u(t,-h_0)=\tilde u(t,\tilde h(t))=0,\ \ &0<t\leq T,\\[1mm]
 \tilde h'(t)\leq-\mu\tilde u_x(t,\tilde h(t)),\ \ &0<t\leq T.
 \end{array}\right.
 \eess
Then the solution $(u,v,g,h)$ of $(\ref{1.2})$ satisfies $h\geq\tilde h$ on  $[0,T]$, and $u\geq\tilde u$ on $\overline{\cal D}$.
 \end{lem}

\begin{lem} $($Comparison principle$)$\label{lm3.3}\, 
Let $T>0$, $\tilde g\in C^1([0,T])$ with $\tilde g<0$ on 
 $[0,T]$, $\tilde u\in C(\overline{\cal O})\cap C^{1,2}({\cal O})$ with
${\cal O}=\{0<t\leq T,\ \tilde g(t)<x<h_0\}$.
Suppose that $(\tilde u, \tilde g)$ satisfies $\tilde g(0)\geq -h_0$, 
$u_0(x)\geq\tilde u(0,x)$ on $[\tilde g(0),h_0]$, and
 \bess
 \left\{\begin{array}{ll}
  \tilde u_t-\tilde u_{xx}\leq\tilde u(1-\tilde u),\ \ &0<t\leq T,\, \ \tilde
g(t)<x<h_0,\\[1mm]
  \tilde u(t,\tilde g(t))=\tilde u(t,h_0)=0,\ \ &0<t\leq T,\\[1mm]
 \tilde g'(t)\geq-\mu\tilde u_x(t,\tilde g(t)),\ \ &0<t\leq T.
 \end{array}\right.
 \eess
Then the solution $(u,v,g,h)$ of $(\ref{1.2})$ satisfies $g\leq\tilde g$ on $[0,T]$, $u\geq\tilde u$ on $\overline{\cal O}$. 
\end{lem}

\section{Long time behavior of $(u,v)$}
\setcounter{equation}{0}

It follows from Lemma \ref{lm2.1} that $x=g(t)$ is monotonic decreasing and $x=h(t)$ is
monotonic increasing. Therefore, there exist $g_\infty\in[-\infty,0)$ and $h_\infty\in(0,\infty]$ such that $\lim_{t\to\infty} g(t)=g_\infty$ and  $\lim_{t\to\infty} h(t)=h_\infty$.
To discuss the long time behavior of $(u,v)$, we first derive an estimate.

\begin{theo}\label{th4.1} \ Let $(u,v,g,h)$ be the solution of $(\ref{1.2})$. If
$h_\infty-g_\infty<\infty$, then there exists a constant $K>0$,
such that
 \bes
 \|u(t,\cdot)\|_{C^1[g(t),h(t)]}\leq K,\ \forall \ t>1,\label{4.1}\\[1mm]
 \lim\limits_{t\to\infty}g'(t)=\lim\limits_{t\to\infty}h'(t)=0.\label{4.2}
 \ees
\end{theo}

{\bf Proof.}  Introduce new functions $w(t,y)$ and $z(t,y)$ by
 \[w(t,y)=u\kk(t,\frac{(h(t)-g(t))y+h(t)+g(t)}{2}\rr),\ \ \
 z(t,y)=v\kk(t,\frac{(h(t)-g(t))y+h(t)+g(t)}{2}\rr).\]
Clearly, $w$ and $z$ satisfy the following initial boundary value problem in an
interval $-1\leq y\leq 1$ with fixed boundary $y=\pm1$:
\be
\left\{\begin{array}{ll}
w_t=\varphi(t)w_{yy}+\psi(t,y)w_y+w(1-w+az),\ \ &t>0,\ \ |y|<1,\\[1mm]
w(t,-1)=w(t,1)=0, &t>0,\\[1mm]
w(0,y)=u(h_0y),&|y|\leq 1,
\end{array}\right.\label{4.3}
\ee
where
 $$\varphi(t)=\frac{4}{(h(t)-g(t))^2},\qquad
\psi(t,y)=\frac{(h'(t)-g'(t))y+h'(t)+g'(t)}{h(t)-g(t)}.$$
By Proposition \ref{pa.1}, we have
$\|w\|_{C^{\frac{1+\alpha}{2},1+\alpha}([1,\infty)\times[-1,1])}<K_0$ for some positive constant $K_0$. Remember $u_x(t,x)=\frac{2}{h(t)-g(t)}w_y(t,y)$. There exists a positive constant $K$ such that
  $$\|u(t,\cdot)\|_{C^1([g(t),\,h(t)])}<K,\ \ \forall \ t\geq 1.$$

We next prove $\lim_{t\to\infty}g'(t)=0$.
Note that $\|w_y(\cdot,-1)\|_{C^{\alpha\over 2}([1,\infty))}<K_0$,
$-M_3<g'(t)<0$ and
$g'(t)=-\mu u_x(t,g(t))=-\frac{2\mu}{h(t)-g(t)}w_y(t,-1)$, it yields
$\|g'\|_{C^{\alpha\over 2}([1,\infty))}<L$,
where $L$ depends on $K_0$ and $M_3$. In view of $g'(t)<0$ and $g_\infty>-\infty$, it is easily to derive that $\lim_{t\to\infty}g'(t)=0$. Analogously, we can obtain
$\lim_{t\to\infty}h'(t)=0$. The proof is complete. \ \ \ \ \fbox{}

\subsection{Vanishing case {\rm(}$h_\infty-g_\infty<\infty${\rm)}}

\begin{theo}\label{th4.2} \ Let $(u,v,g,h)$ be any solution of $(\ref{1.2})$. If
$h_\infty-g_\infty<\infty$, then
 \bes
 &\dd\lim_{t\to\infty}\|u(t,\cdot)\|_{C([g(t),h(t)])}=0,&\label{4.4}\\[1mm]
 &\dd\lim_{t\to\infty}v(t,x)=b\ \ \mbox{uniformly\, on\, the\, compact\, subset\, of } \, \mathbb{R}.&\label{4.5}\ees
This result shows that if the predator can not spread into the whole space, then it
will die out eventually.
\end{theo}

From the results of Section $5$ we shall see that the reason leading to the predator species disappears eventually are three aspects: (a) the initial habitat $[-h_0, h_0]$ of the predator  is too narrow, (b) the initial data $u_0(x)$ of the predator is too small, or (c) the moving parameter/coefficient $\mu$ of free boundaries is too small.

{\bf Proof of Theorem \ref{th4.2}}

{\it Step 1: Proof of $(\ref{4.4})$}. On the contrary we assume that
there exist $\varepsilon>0$ and $\kk\{(t_j,x_j)\rr\}_{j=1}^{\infty}$, with
$g(t_j)<x_j<h(t_j)$ and  $t_j\to\infty$ as $j\to\infty$, such that
 \bes u(t_j,x_j)\geq 3\varepsilon,\ j=1,2,\cdots.\label{4.6}\ees
Since $g_\infty<x_j<h_\infty$, there are a subsequence of $\{x_j\}$, noted by itself,
and $x_0\in[g_\infty,h_\infty]$, such that $x_j\to x_0$ as $j\to\infty$.
We claim that $x_0\in(g_\infty,h_\infty)$. If $x_0=g_\infty$, then $x_j-g(t_j)\to 0$ as $j\to\infty$. By use of the inequality (\ref{4.6}) firstly and the inequality
(\ref{4.1}) secondly, it is deduced that
 \bess
 \dd\frac{4\varepsilon}{x_j-g(t_j)}\leq\frac{u(t_j,x_j)}{x_j-g(t_j)}=
 \frac{u(t_j,x_j)-u(t_j,g(t_j))}{x_j-g(t_j)}=u_x(t_j,\bar x_j)\leq K,\eess
where $\bar x_j\in[g(t_j),x_j]$. It is a contradiction as $x_j-g(t_j)\to 0$.
Similarly, we can ensure $x_0<h_\infty$.

By use of (\ref{4.1}) and (\ref{4.6}), there exists $\delta>0$ such that  $[x_0-\delta , \,
x_0+\delta]
\subset(g_\infty,h_\infty)$ and
 $$u(t_j,x)\geq 2\varepsilon,\ \ \ \forall\ x\in[x_0-\delta, \, x_0+\delta]$$
for all large $j$. As $g(t_j)\to g_\infty$ and $h(t_j)\to h_\infty$ as
$j\to\infty$, without loss of generality we may think that $g(t_j)<x_0-\delta$ and
$h(t_j)>x_0+\delta$ for all $j$.

Let $l_j(t)=x_0-\delta-(t-t_j)$, $r_j(t)=x_0+\delta+t-t_j$. Then $l_j(t_j)>g(t_j)$
and $r_j(t_j)<h(t_j)$. Set
 $$\tau_j=\inf\kk\{t>t_j:\, g(t)=l_j(t),\ or\ h(t)=r_j(t)\rr\}.$$
Since $h_\infty<\infty$, $g_\infty>-\infty$, and $l_j(t)\to-\infty$ and
$r_j(t)\to\infty$ as $t\to\infty$, we see that $t^*_j<\infty$. It is easy to
obtain that $\tau_j<t_j-\delta+(h_\infty-g_\infty)/2$. Without loss of
generality, we assume that $h(\tau_j)=r_j(\tau_j)$ for all $j$. This implies
 \bes
 g(t)\leq l_j(t)<r_j(t)\leq h(t)\ \ \ \mbox{in}\ \  [t_j,\tau_j].
 \label{4.7}\ees
Define $y_j(t,x)=(\pi-\theta)\dd\frac{x-x_0}{\delta+t-t_j}$ and
 $$ u_j(t,x)=\varepsilon {\rm e}^{-k(t-t_j)}[\cos y_j(t,x)+\cos
\theta],\ \ \ (t,x)\in\overline\Omega_j,$$
where $\theta\ (\theta<\pi/8)$ and $k$ are positive constants to be chosen
later, and
   $$\Omega_j=\{(t,x):t_j<t<\tau_j, \ l_j(t)<x<r_j(t)\}.$$
It is obvious that $ u_j(t,l_j(t))=0= u_j(t,r_j(t))$, and
$|y_j(t,x)|\leq \pi-\theta$ for $(t,x)\in\overline\Omega_j$, the latter
implies $ u_j(t,x)\geq 0$ in $\Omega_j$.

We want to compare $u(t,x)$ and $ u_j(t,x)$ in $\overline\Omega_j$.
Thanks to (\ref{4.7}), it follows that
 $$u(t,l_j(t))\geq 0= u_j(t,l_j(t)), \ \ u(t,r_j(t))\geq 0=
u_j(t,r_j(t))\ \ \ \mbox{for} \ \ t\in[t_j,\tau_j].$$
On the other hand, it is obvious that
 $$u(t_j,x)\geq 2\varepsilon\geq u_j(t_j,x)\ \ \
  \mbox{for} \ \ x\in[x_0-\delta,x_0+\delta].$$
Thus, if the positive constants $\theta$ and $k$ can be chosen independent of
$j$ such that
 \bes u_{jt}- u_{jxx}- u_j(1- u_j)\leq 0
  \ \ \ \mbox{in} \ \ \Omega_j,
  \label{4.8}\ees
it can be deduced that $ u_j(t,x)\leq u(t,x)$ for $(t,x)\in\Omega_j$ by
applying the maximum principle to $u- u_j$ over $\Omega_j$.
Since $u(\tau_j,h(\tau_j))=0= u_j(\tau_j,r_j(\tau_j))$ and
$h(\tau_j)=r_j(\tau_j)$, it follows that $u_x(\tau_j,h(\tau_j))\leq
u_{jx}(\tau_j,r_j(\tau_j))$. Thanks to $\ep<\pi/8$ and $\delta+\tau_j-t_j<(h_\infty-g_\infty)/2$, we derive
  $$ u_{jx}(\tau_j,r_j(\tau_j))=-\dd\frac{\varepsilon(\pi-\theta)}
{\delta+\tau_j-t_j}{\rm e}^{-k(\tau_j-t_j)}\sin(\pi-\theta)
\leq-\dd\frac{7\varepsilon\pi}{4(h_\infty-g_\infty)}{\rm e}^{-k(h_\infty-g_\infty)}
 \sin \theta.$$
Note the boundary condition $-\mu u_x(\tau_j,h(\tau_j))=h'(\tau_j)$, one has immediately
  $$h'(\tau_j)\geq\dd\frac{7\mu\varepsilon\pi}{4(h_\infty-g_\infty)}
  {\rm e}^{-k(h_\infty-g_\infty)}\sin \theta,$$
which implies $\limsup_{t\to\infty}|h'(t)|>0$ since $\lim_{j\to\infty}\tau_j\to\infty$. This contradicts to (\ref{4.2}), and (\ref{4.4}) is obtained.

\vskip 4pt We claim that (\ref{4.8}) holds so long as $\theta$ and $k$ satisfy
\bes
 \theta<\frac{\pi}{8}, \ \ \sin\theta<
 \frac{3\delta^2\pi}{(h_\infty-g_\infty)^3}, \ \ 
 k>\frac{\pi(h_\infty-g_\infty)}{2\delta^2(\cos\theta-\cos
2\theta)}+2\varepsilon +\kk(\frac{\pi}{\delta}\rr)^2.\label{4.10}
  \ees
In fact, a series of computations indicate that, for $(t,x)\in\Omega_j$,
  $$\begin{array}{ll}
 & u_{jt}- u_{jxx}- u_j(1-u_j)\\[1mm]
 =&-k u_j-\varepsilon {\rm e}^{-k(t-t_j)}y_{jt}\sin
y_j+\varepsilon {\rm e}^{-k(t-t_j)}y_{jx}^2\cos y_j-
u_j(1- u_j)\\[1mm]
\leq&\dd\kk(2\varepsilon +y_{jx}^2-k\rr)
u_j-\varepsilon {\rm e}^{-k(t-t_j)}y_{jx}^2\cos
\theta-\varepsilon {\rm e}^{-k(t-t_j)}y_{jt}\sin y_j\\[1mm]
\leq&\kk(2\varepsilon +(\pi/\delta)^2-k\rr) u_j
-4\varepsilon {\rm e}^{-k(t-t_j)}\kk(\frac{\pi-\theta}
{h_\infty-g_\infty}\rr)^2\cos \theta
+\frac{\varepsilon\pi|x-x_0|}{\delta^2}{\rm e}^{-k(t-t_j)}|\sin y_j|\\[1mm]
:=&I(t,x).
  \end{array}$$
Obviously, $2\varepsilon +(\pi/\delta)^2-k<0$ by the third inequality of (\ref{4.10}). Since $|y_j(t,x)|\leq \pi-\theta$ in $\overline\Omega_j$, we can decompose $\Omega_j=D_j\bigcup E_j$ with
\bess
  D_j&=&\kk\{(t,x)\in\Omega_j:\, t_j<t<\tau_j,\
 \pi-2\theta<|y_j(t,x)|<\pi-\theta\rr\},\\[1mm]
 E_j&=&\kk\{(t,x)\in\Omega_j:\, t_j<t<\tau_j,
 \ |y_j(t,x)|\leq\pi-2\theta\rr\}.\eess
It is obvious that $|\sin y_j(t,x)|\leq \sin 2\theta$ in $D_j$,
$\cos y_j(t,x)\geq-\cos 2\theta$ in $E_j$. Note that
$ u_j(t,x)\geq 0$ and $|x-x_0|\leq (h_\infty-g_\infty)/2$ in $\Omega_j$,
in view of (\ref{4.10}), we conclude
 $$ I(t,x)\leq \varepsilon {\rm e}^{-k(t-t_j)}\kk(
 -\dd\frac{3 \pi^2}{(h_\infty-g_\infty)^2}\cos \theta
 +\dd\frac{\pi(h_\infty-g_\infty)}{2\delta^2}\sin 2\theta\rr)<0$$
when $(t,x)\in D_j$, and
 $$I(t,x)\leq \varepsilon {\rm e}^{-k(t-t_j)}\kk(\big(2\varepsilon
 +({\pi}/{\delta})^2-k\big)(\cos\theta-\cos 2\theta)
 +\dd\frac{\pi(h_\infty-g_\infty)}{2\delta^2}\rr)<0$$
when $(t,x)\in E_j$. Therefore, (\ref{4.8}) holds.

\vskip 2pt {\it Step 2: Proof of (\ref{4.5})}. By the comparison principle, $v(t,x)\leq\bar v(t)$ for all $t\in[0,\infty)$ and $x\in\mathbb{R}$, where
 $$\bar v(t)=b{\rm e}^{bt}\kk({\rm e}^{bt}-1
 +\dd\frac{b}{\|v_0\|_\infty}\rr)^{-1},$$
which is the solution of the ODE problem
  $$\bar v'(t)=\bar v(b-\bar v),\ \ t>0;\ \ \bar v(0)=\|v_0\|_\infty.$$
Since $\lim_{t\to\infty}\bar v(t)=b$, it is deduced that
$\limsup_{t\to\infty}v(t,x)\leq\lim_{t\to\infty}\bar v(t)=b$
uniformly for $x\in\mathbb{R}$.

On the other hand, note that (\ref{4.4}) and $u(t,x)\equiv 0$ for $t>0$, $ x\not\in(g(t), h(t))$, we see that for any given $0<\sigma\ll 1$, there exists $T_\sigma>0$ such that $u(t,x)<\sigma$ for $t>T_\sigma$ and $x\in\mathbb{R}$.
For any given $\ep>0$ and $L>0$, let $l_\ep$ be determined by Proposition \ref{pB.1}. Then $v$ satisfies
  $$ \left\{\begin{array}{ll}
  v_t-D v_{xx}\geq v(b-v-c\sigma), \ \ &t>T_\sigma, \ \ -l_\ep<x<l_\ep,\\[1mm]
 v(t,\pm l_\ep)>0,&t\ge T_\sigma,\\[1mm]
 v(T_\sigma,x)>0,&-l_\ep\leq x\leq l_\ep.
 \end{array}\right.$$
Thanks to Proposition \ref{pB.1}, $\liminf_{t\to\infty}v(t,x)\geq (b-c\sigma)-\ep$ uniformly on $[-L,L]$. By the arbitrariness of $\ep$ and $L$, we derive that $\liminf_{t\to\infty}v(t,x)\geq b-c\sigma$ uniformly in the compact subset of $\mathbb{R}$.
Since $\sigma>0$ is arbitrary, it follows that
$\liminf_{t\to\infty} v(t,x)\geq b$ uniformly in any bounded
subset of $\mathbb{R}$. The proof is complete.  \ \ \ \ \fbox{}

\subsection{Spreading case {\rm(}$h_\infty-g_\infty=\infty${\rm)}}

We first provide a proposition which asserts the equivalence of $h_\infty-g_\infty=\infty$ and $g_\infty=-\infty$, $h_\infty=\infty$.

\begin{prop}\label{p4.1} \ If $h_\infty-g_\infty=\infty$, then $g_\infty=-\infty$ and $h_\infty=\infty$.
\end{prop}

{\bf Proof.}  Since $h_\infty-g_\infty=\infty$, there exists $T>0$ such that
$h(T)-g(T)>\pi$. Choose a function $\tilde u_0(x)$ satisfying $\tilde
u_0\in C^2[g(T),h(T)]$, $\tilde u_0(x)\leq u_0(x,T)$ in  $[g(T),h(T)]$, $\tilde
u_0(x)>0$ in $(g(T),h(T))$ and $\tilde u_0(g(T))=\tilde u_0(h(T))=0$. Consider the
following problem
  $$ \left\{\begin{array}{ll}
  \tilde u_t-\tilde u_{xx}=\tilde u(1-\tilde u),\ \ &t>T,\ \ g(T)<x<\tilde
h(t),\\[1mm]
 \tilde u(t, g(T))=0=\tilde u(t,\tilde h(t)),\ \ &t\ge T,\\[1mm]
 \tilde h'(t)=-\mu\tilde u_x(t,\tilde h(t)),\ \ &t\ge T,\\[1mm]
 \tilde h(T)=h(T), \ \tilde u(T,x)=\tilde u_0(x), \ &g(T)\leq x\leq h(T).
 \end{array}\right.$$
By Theorem  2.7 of \cite{KY}, this problem has a unique solution $(\tilde u, \tilde
h)$ and exists for all $t\geq T$, and by Theorem 4.2 of \cite{KY},
$\lim_{t\to\infty}\tilde h(t)=\infty$. In view of Lemma \ref{lm3.2}, it concludes
$h(\infty)=\lim_{t\to\infty} h(t)\geq\lim_{t\to\infty}\tilde h(t)=\infty$. 
Similarly, by Lemma \ref{lm3.3} we can obtain $g(\infty)=-\infty$. \ \ \ \ \fbox{}

\vskip 4pt We deal with the weakly hunting case $b>c$ and $ac<1$ firstly.

\begin{theo}\label{th4.3} Assume that $g_\infty=-\infty$, $h_\infty=\infty$. For the
weakly hunting case $b>c$ and $ac<1$, we have
 \bes
 \lim_{t\to\infty}u(t,x)=\frac{1+ab}{1+ac},\ \ \ \
\lim_{t\to\infty}v(t,x)=\frac{b-c}{1+ac}
 \label{4.11}
 \ees
uniformly in any compact subset of $\mathbb{R}$.
 \end{theo}

{\bf Proof.} For any given $L>0$ and $0<\ep\ll 1$, let $l_\ep$ be given by Proposition \ref{pB.1} with $d=1$, $\beta=1$ and $\theta=1$.
In view of $g_\infty=-\infty$ and $h_\infty=\infty$, there exists $T_0>0$ such that
 \[g(t)<-l_\ep, \ \ h(t)>l_\ep, \ \ \
  \forall \ t\geq T_0.\]
Notice that $v>0$, we see that $u$ satisfies
 \bess
 \left\{\begin{array}{lll}
  u_t-u_{xx}> u(1-u),\ &t\geq T_0, \ \ x\in[-l_\ep,l_\ep],\\[1mm]
   u(t,\pm l_\ep)>0,\ \ &t\geq T_0.
 \end{array}\right.
 \eess
Since $u(T_0,x)>0$ in $[-l_\ep,l_\ep]$, applying Proposition \ref{pB.1}, it arrives at
 \[\liminf_{t\to\infty}u(t,x)\geq 1-\ep\ \ \ \mbox{uniformly\, on }\, [-L,L].\]
By the arbitrariness of $\ep$ and $L$,
  \bes\liminf_{t\to\infty}u(t,x)\geq 1:=\underline{u}_1\ \ \mbox{uniformly\, on\, the\, compact\,  subset\, of } \, \mathbb{R}. \label{4.12}\ees

Let $M=\max\{M_1,M_2\}$, where $M_i$ is
determined by Lemma \ref{lm2.1}, $i=1,2$. For any given $L>0$, $0<\delta\ll
1$ and $0<\ep\ll 1$, let $l_\ep$ be given by Proposition \ref{pB.2} with $d=D$,
$\beta=b-c\kk(\underline{u}_1-\delta\rr)$, $\theta=1$ and $k=M$. In view of (\ref{4.12}), there exists $T_1>0$ such that $u(t,x)\geq \underline{u}_1-\delta$ for all
$t\geq T_1$ and $x\in[-l_\ep,l_\ep]$. Therefore, $v$ satisfies
 \bess
 \left\{\begin{array}{lll}
  v_t-Dv_{xx}\leq v\kk[b-v-c\kk(\underline{u}_1-\delta\rr)\rr],
  \ &t\geq T_1, \ \ x\in[-l_\ep,l_\ep],\\[1mm]
  v(t,\pm l_\ep)\leq M,\ \ &t\geq T_1.
 \end{array}\right.
 \eess
As $v(T_1,x)>0$ in $[-l_\ep,l_\ep]$, in view of Proposition \ref{pB.2}, it yields
 \[\limsup_{t\to\infty}v(t,x)\leq b-c\big(\underline{u}_1
 -\delta\big)-\ep\ \ \ \mbox{uniformly\, on }\, [-L,L].\]
The arbitrariness of $\ep$, $L$ and $\delta$ imply that
  \bes\limsup_{t\to\infty}v(t,x)\leq b-c\underline{u}_1:=\bar v_1\ \ \mbox{uniformly\, on\, the\, compact\,  subset\, of } \, \mathbb{R}. \label{4.13}\ees

For any given $L>0$, $0<\delta\ll 1$ and $0<\ep\ll 1$, let $l_\ep$ be given by
Proposition \ref{pB.2} with $d=1$, $\beta=1+a\kk(\bar v_1+\delta\rr)$, $\theta=1$ and
$k=M$. Taking into account (\ref{4.13}), and $g_\infty=-\infty$ and $h_\infty=\infty$, there is $T_2>0$ such that
 \[v(t,x)\leq \bar v_1+\delta,\ \ g(t)<-l_\ep, \ \ h(t)>l_\ep, \ \ \
  \forall \ t\geq T_2, \ \ x\in[-l_\ep,l_\ep].\]
Hence, $u$ satisfies
 \bess
 \left\{\begin{array}{lll}
  u_t-u_{xx}\leq u\kk[1-u+a\kk(\bar v_1+\delta\rr)\rr],\ &t\geq T_2,
   \ \ x\in[-l_\ep,l_\ep],\\[1mm]
   u(t,\pm l_\ep)\leq M,\ \ &t\geq T_2.
 \end{array}\right.
 \eess
By the same argument as above, one gets
\bes\limsup_{t\to\infty}u(t,x)\leq 1+a\bar v_1:=\bar u_1\ \
\mbox{uniformly\, on\, the\, compact\,  subset\, of } \, \mathbb{R}. \label{4.14}\ees

For any given $L>0$, $0<\delta\ll 1$ and $0<\ep\ll 1$, let $l_\ep$ be given by Proposition \ref{pB.1} with $d=D$, $\beta=b-c\kk(\bar u_1+\delta\rr)$ and $\theta=1$. According to (\ref{4.14}), there is $T_3>0$ such that $u(t,x)\leq\bar u_1+\delta$ for
all $t\geq T_3$ and $x\in[-l_\ep,l_\ep]$. Hence, $v$ satisfies
 \bess
 \left\{\begin{array}{lll}
  v_t-Dv_{xx}\geq v\kk[b-v-c\kk(\bar u_1+\delta\rr)\rr],\ &t\geq T_3,\ \ x\in[-l_\ep,l_\ep],\\[1mm]
  v(t,\pm l_\ep)\geq 0,\ \ &t\geq T_3.
 \end{array}\right.
 \eess
Similar to the above,
 \bes\liminf_{t\to\infty}v(t,x)\geq b-c\bar{u}_1:=\underline{v}_1\ \ \mbox{uniformly\, on\, the\,
compact\,  subset\, of } \, \mathbb{R}. \label{4.15}\ees

For any given $L>0$, $0<\delta\ll 1$ and $0<\ep\ll 1$, let $l_\ep$ be given by
Proposition \ref{pB.1} with $d=1$, $\beta=1+a\kk(\underline{v}_1-\delta\rr)$
and $\theta=1$. By virtue of (\ref{4.15}), and $g_\infty=-\infty$ and $h_\infty=\infty$, there is $T_4>0$ such that
 \[v(t,x)\geq \underline{v}_1-\delta,\ \ g(t)<-l_\ep, \ \ h(t)>l_\ep,\ \ \
  \forall \ t\geq T_4, \ \ x\in[-l_\ep,l_\ep].\]
Thus, $u$ satisfies
 \bess
 \left\{\begin{array}{lll}
  u_t-u_{xx}\geq u\kk[1-u+a\kk(\underline{v}_1-\delta\rr)\rr],
  \ &t\geq T_4,\ \ x\in[-l_\ep,l_\ep],\\[1mm]
   u(t,\pm l_\ep)\geq 0,\ \ &t\geq T_4.
 \end{array}\right.
 \eess
Same as above,
\bess\liminf_{t\to\infty}u(t,x)\geq
(1+a\underline{v}_1):=\underline{u}_2\ \ \mbox{uniformly\, on\, the\,
compact\,  subset\, of } \, \mathbb{R}.\eess

Repeating the above procedure, we can find four sequences $\{\underline{u}_i\}$,
$\{\underline{v}_i\}$, $\{\bar{u}_i\}$ and $\{\bar{v}_i\}$, such that, for all $i$,
 \bes
 \underline{u}_i\leq\liminf_{t\to\infty}u(t,x)\leq
 \limsup_{t\to\infty}u(t,x)\leq\bar{u}_i,\quad
 \underline{v}_i\leq\liminf_{t\to\infty}v(t,x)\leq
\limsup_{t\to\infty}v(t,x)\leq\bar{v}_i
\label{4.16}\ees
uniformly in the compact subset of $\mathbb{R}$. Moreover, these sequences can be
determined by the following iterative formulas:
 \bes
 \underline{u}_1=1, \ \ \bar v_i=b-c\underline{u}_i, \ \
\bar u_i=1+a\bar v_i,\ \
 \underline{v}_i=b-c\bar{u}_i, \ \
\underline{u}_{i+1}=1+a\underline{v}_i,\ \ i=1,2,\cdots.
\label{4.17}\ees

Denote $A=b-c$ and $q=ac$, then $A>0$, $0<q<1$. By the direct calculation,
\bess
\bar v_1=A, \ \ \bar u_1=1+aA, \ \
\underline{v}_1=A(1-q), \ \  \underline{u}_2=\dd 1+aA(1-q), \ \
\underline{v}_2=A(1-q+q^2).\eess
Using the inductive method we have the following expressions:
 \bess
 \bar v_i=A\kk(1-q+q^2-\cdots+q^{2i-4}-q^{2i-3}+q^{2i-2}\rr),
 \ \  \underline{v}_i=\bar v_i-Aq^{2i-1},
 \ \ \ i\geq 3.\eess
Because $0<q<1$, one has
 \bes\lim_{i\to\infty}\bar v_i=\lim_{i\to\infty}\underline{v}_i=\frac A{1+q}
 =\frac{b-c}{1+ac}.\label{4.18}\ees
This fact combines with (\ref{4.17}) yields that
  \bes\lim_{i\to\infty}\bar u_i=\lim_{i\to\infty}\underline{u}_i
 =\frac{1+ab}{1+ac}.\label{4.19}\ees
The limits (\ref{4.11}) are followed from (\ref{4.16}), (\ref{4.18}) and  (\ref{4.19}). \ \ \ \ \fbox{}

\vskip 4pt For the strongly hunting case: $b\leq c$, similar to the above, the following theorem can be obtianed.

\begin{theo}\label{th4.4} Assume that $g_\infty=-\infty$ and $h_\infty=\infty$. For the
strongly hunting case $b\leq c$, we have
 \[\lim_{t\to\infty}u(t,x)=1,\ \ \ \ \lim_{t\to\infty}v(t,x)=0\]
uniformly in any compact subset of $\mathbb{R}$.
 \end{theo}

\section{The criteria governing spreading and vanishing}
\setcounter{equation}{0}

We first give a necessary condition for vanishing.

\begin{theo}\label{th5.1} \ Let $(u,v,g,h)$ be any solution of {\rm(\ref{1.2})}. If
$h_\infty-g_\infty<\infty$,
then
  \bes
 h_\infty-g_\infty\leq\pi\sqrt{1/(1+ab)}:=\Lambda.\label{5.1}\ees
Hence, $h_0\geq \Lambda/2$ implies
$h_\infty-g_\infty=\infty$ due to $g'(t)<0$ and $h'(t)>0$ for $t>0$.
\end{theo}

{\bf Proof.} By Theorem \ref{th4.2}, if $h_\infty-g_\infty<\infty$ then
$\lim_{t\to\infty}\|u(t,\cdot)\|_{C[g(t),h(t)]}=0$ and $\lim_{t\to\infty}v(t,x)=b$ uniformly in the bounded subset of $\mathbb{R}$. We assume $h_\infty-g_\infty>\Lambda$ to get a contradiction. For any small $\varepsilon>0$, there exists $\tau\gg 1$ such that
  \bess
  &v(t,x)\dd\geq b-\varepsilon/a:=A_\varepsilon,\ \ \
  \forall \ t\geq\tau, \ x\in[g_\infty,h_\infty],&\\[1mm]
 &h(\tau)-g(\tau)>\max\big\{2h_0, \
  \pi\sqrt{1/(1+ab-\varepsilon)}\big\}.&\eess
Set $l_1=g(\tau)$ and $l_2=h(\tau)$, then
$l_2-l_1>\pi\sqrt{1/(1+ab-\varepsilon)}$.
Let $w$ be the positive solution of the following initial boundary value
problem with fixed boundary:
  $$\left\{\begin{array}{ll}
 w_t=w_{xx}+w\dd\kk(1-w+aA_\varepsilon\rr),
  \ \ &t>\tau, \ \, l_1<x<l_2,\\[1mm]
 w(t,l_1)=w(t,l_2)=0,&t>\tau,\\[1mm]
  w(\tau,x)=u(\tau,x),&l_1<x<l_2.
  \end{array}\right.$$
By the comparison principle,
 $$w(t,x)\leq u(t,x)\ \ \ \mbox{for} \ \ t\geq\tau,\ l_1\leq x\leq l_2.$$
Since $1+aA_\varepsilon>\kk[{\pi}/(l_2-l_1)\rr]^2$,
it is well known that $w(t,x)\to \theta(x)$ as $t\to\infty$ uniformly in the
compact subset of $(l_1,l_2)$, where $\theta$ is the unique positive solution of
  \[\left\{\begin{array}{ll}
 \theta_{xx}+\theta\dd\kk(1+aA_\varepsilon-\theta\rr)=0,\ \ &l_1<x<l_2,\\[1mm]
 \theta(l_1)=\theta(l_2)=0.&\end{array}\right.\]
Hence, $\liminf_{t\to\infty} u(t,x)\geq\lim_{t\to\infty}w(t,x)=\theta(x)>0$ in
$(l_1,l_2)$. This is a contradiction to (\ref{4.4}). Consequently, (\ref{5.1}) holds.
\ \ \ \ \fbox{}

\vskip 2pt By Theorem \ref{th5.1} and Proposition \ref{p4.1}, $h_0\geq\Lambda/2$
implies $g_\infty=-\infty$ and $h_\infty=\infty$.

\vskip 2pt Now we discuss the case $h_0<\Lambda/2$.

\begin{lem}\label{lmq5.1} \ Suppose that $h_0<\Lambda/2$. If
  \bess
  \mu\geq\mu^0:=\max\kk\{1,\|u_0\|_\infty\rr\}\kk(\pi^2
  -4h^2_0\rr)\kk(2\dd\int_{-h_0}^{h_0}(x+h_0)u_0(x)dx\rr)^{-1},
  \eess
then $g_\infty=-\infty$ and $h_\infty=\infty$.
\end{lem}

{\bf Proof.} Consider the following auxiliary problem
$$
 \left\{\begin{array}{ll}
  \underline u_t-\underline u_{xx}=\underline u(1-\underline u),
  \ &t>0, \ \ -h_0<x<\underline h(t),\\[1mm]
 \underline u(t,-h_0)=0,\ \ \underline u(t,\underline h(t))=0,\ &t\ge 0,\\[1mm]
 \underline h'(t)=-\mu\underline u_x(t,\underline h(t)), \ &t\ge 0,\\[1mm]
 \underline h(0)=h_0, \ \underline u(0,x)=u_0(x),\ \ &-h_0\leq x\leq h_0.
 \end{array}\right.$$
It follows from Lemma \ref{lm3.2} that
  $$\underline h(t)\leq h(t),\ \ \underline u(t,x)\leq u(t,x)
  \ \ {\rm for}\ t>0\ {\rm and}\ -h_0<x<\underline h(t).$$
Recall that $2h_0<\Lambda<\pi$ and $\mu\geq\mu^0$, by Proposition 4.8 of \cite{KY}, we have $\underline
h(\infty)=\infty$. Therefore, $h_\infty=\infty$. Similarly, $g_\infty=-\infty$. The
proof is finished. \ \ \ \ \fbox{}

\begin{lem}\label{lmq5.2} \ Assume that $h_0<\Lambda/2$. Then
there exists $\mu_0>0$, depending also on $u_0(x)$ and $v_0(x)$, such that
$h_\infty-g_\infty<\infty$ when $\mu\leq\mu_0$.
\end{lem}

{\bf Proof.} We are going to construct a suitable supper solution to (\ref{1.2})
and then apply Lemma \ref{lm3.1}.
Obviously, the function
  $$\bar v(t)=\dd b{\rm e}^{bt}\kk({\rm e}^{bt}-1
 +\dd\frac{b}{\|v_0\|_\infty}\rr)^{-1}$$
satisfies
 $$\left\{\begin{array}{ll}
\bar v_t-D\bar v_{xx}\geq \bar v(b-\bar v),\ \ &x\in\mathbb{R},\ t>0,\\[1mm]
\bar v(0,x)\geq v_0(x),&x\in\mathbb{R}.
\end{array}\right.$$
Denote $\vartheta=\frac 12 h_0+\frac 14\Lambda$, then
$h_0<\vartheta<\Lambda/2$. Inspired by \cite{RT}, we define
 \bess
 &f(t)=M\exp\dd\kk\{\int_0^t\kk[1+a\bar
v(s)-\kk(\dd\frac{\pi}{2\vartheta}\rr)^2\rr]ds\rr\},&\\[1mm]
 &\eta(t)=\dd\kk(h_0^2(1+\delta)^2+\dd\mu\pi\int_0^tf(s)ds\rr)^{1/2},\ \ t\geq 0;
  \ \ \ w(y)=\cos\dd\frac{\pi y}{2},\ \ -1\leq y\leq 1,&\\[1mm]
&\bar u(t,x)=f(t) w\kk(\dd\frac{x}{\eta(t)}\rr),\ \ \
    t\geq 0,\ \ -\eta(t)\leq x\leq \eta(t),&\eess
where $\delta\ll 1$ is a fixed positive constant such that $\vartheta>h_0(1+\delta)$ and $M$
is a positive constant to be chosen later. Clearly, $\eta'(t)>0$ for $t\geq 0$ and
  \bes\frac{f'(t)}{f(t)}=1+a\bar v(t)-\kk(\dd\frac{\pi}
  {2\vartheta}\rr)^2\ \ \ {\rm for}\ t>0.\label{5.2}\ees
Remember that $\vartheta<\Lambda/2$ and $\lim_{t\to\infty}\bar v(t)=b$, we have
$1+a\bar v(t)-\kk(\frac{\pi}{2\vartheta}\rr)^2<0$ for $t$ large enough, and then
$\int_0^tf(s)ds$ is uniformly bounded in $[0,\infty)$.

Let
  $$\mu_0=\frac{\vartheta^2-h_0^2(1+\delta)^2}{\pi\int_0^\infty f(t)dt}.$$
When $0<\mu\leq\mu_0$, it is obvious that $\vartheta\geq \eta(t)$ for all $t\geq 0$. In view of
(\ref{5.2}), we have that by the direct computation
  $$\begin{array}{rl}
\bar u_t-\bar u_{xx}-\bar u(1-\bar u+a\bar v)
=&f'w-fw'\frac{x\eta'}{\eta^2}+f\big(\frac{\pi}{2\eta}\big)^2w
-fw(1-fw+a\bar v)\\[2mm]
\geq&fw\Big[\frac{f'}{f}+\big(\frac{\pi}{2\eta}\big)^2-1-a\bar
v\Big]\\[2mm]
 =&\frac{\pi^2}4fw\big(\eta^{-2}-\vartheta^{-2}\big)\geq 0
  \end{array}$$
for all $t>0$ and $-\eta(t)<x<\eta(t)$. On the other hand,
 $$\eta'(t)=\dd\frac{\mu\pi}{2\eta(t)}f(t), \ \
 \bar u_x(t,-\eta(t))=\dd\frac{\pi}{2\eta(t)}f(t), \ \
 \bar u_x(t,\eta(t))=-\dd\frac{\pi}{2\eta(t)}f(t),$$
which imply
 $$-\eta'(t)=-\mu\bar u_x(t,-\eta(t)), \ \ \ \eta'(t)=-\mu\bar u_x(t,\eta(t)).$$

Choose $M$ is so large that $u_0(x)\leq M\cos\frac{\pi x}{2h_0(1+\delta)}$ for
$x\in[-h_0,h_0]$. Then for any $0<\mu\leq\mu_0$, the pair $(\bar u,\bar v)$
satisfies
 $$ \left\{\begin{array}{ll}
  \bar u_t-\bar u_{xx}\geq\bar u(1-\bar u+a\bar v),\ \ &t>0, \ \
 |x|<\eta(t),\\[1mm]
 \bar u(t,\pm\eta(t))=0,\ \ \eta'(t)=\mp\mu \bar u_x(t,\pm\eta(t)), \ &t>0,\\[1mm]
  \bar u(0,x)\geq u_0(x),\ \ &|x|\leq h_0,\\[1mm]
  \eta(0)>h_0.&
  \end{array}\right.$$
Take advantage of Lemma \ref{lm3.1}, $-\eta(t)\leq g(t)$, $\eta(t)\geq
h(t)$ and $u(t,x)\leq\bar u(t,x)$ for $t>0$ and $g(t)\leq x\leq h(t)$. It follows that
  $$g_\infty\geq-\lim_{t\to\infty}\eta(t)>-\vartheta>-\infty,\ \
h_\infty\leq\lim\limits_{t\to\infty}\eta(t)<\vartheta<\infty.$$
The proof is complete. \ \ \ \ \fbox{}

\begin{theo}\label{th5.2} \ Suppose that $h_0<\Lambda/2$. Then there exist $\mu^*\geq\mu_*>0$, depending on $u_0(x)$, $v_0(x)$ and $h_0$, such that $g_\infty=-\infty$ and $h_\infty=\infty$ if $\mu>\mu^*$, and $h_\infty-g_\infty\leq\Lambda$ if $\mu\leq\mu_*$ or $\mu=\mu^*$.
\end{theo}

{\bf Proof.}  The proof is similar to that of Theorem 3.9 of \cite{DLin} and Theorem
4.11 of \cite{KY}. For the convenience to reader we shall give the details. Write $(u_\mu, v_\mu, g_\mu, h_\mu)$ in place of $(u, v, g, h)$ to clarify the
dependence of the solution of (\ref{1.2}) on $\mu$. Define
 $$\Sigma^*=\kk\{\mu>0:\,h_{\mu,\infty}-g_{\mu,\infty}\leq\Lambda\rr\}.$$
By Lemma \ref{lmq5.2} and Theorem \ref{th5.1}, $(0,\mu_0]\subset\Sigma^*$. In view of
Lemma \ref{lmq5.1}, $\Sigma^*\cap[\mu^0,\infty)=\emptyset$. Therefore,
$\mu^*:=\sup\Sigma^*\in[\mu_0,\,\mu^0]$. By this definition and Theorem \ref{th5.1} we
find that $g_{\mu,\infty}=-\infty$ and $h_{\mu,\infty}=\infty$ when $\mu>\mu^*$.
Hence, $\Sigma^*\subset(0,\mu^*]$.

We will show that $\mu^*\in\Sigma^*$. Otherwise, $g_{\mu^*,\infty}=-\infty$ and
$h_{\mu^*,\infty}=\infty$. There exists $T>0$ such that
$h_{\mu^*}(T)-g_{\mu^*}(T)>\Lambda$. Utilizing the continuous
dependence of $(u_\mu, v_\mu, g_\mu, h_\mu)$ on $\mu$, there is $\varepsilon>0$ such
that $h_{\mu}(T)-g_{\mu}(T)>\Lambda$ for
$\mu\in(\mu^*-\varepsilon,\mu^*+\varepsilon)$. It follows that for all such $\mu$,
 \[\lim_{t\to\infty}\big[h_{\mu}(t)-g_{\mu}(t)\big]\geq h_{\mu}(T)-g_{\mu}(T)
 >\Lambda.\]
Therefore, $[\mu^*-\varepsilon,\mu^*+\varepsilon]\cap\Sigma^*=\emptyset$, and
$\sup\Sigma^*\leq\mu^*-\varepsilon$. This contradicts the definition of $\mu^*$.

Define
 $$\Sigma_*=\kk\{\nu:\, \nu\geq\mu_0 \ \mbox{such\ that} \  \,h_{\mu,\infty}-g_{\mu,\infty}\leq\Lambda\ \mbox{for\ all} \ \mu\leq\nu\rr\},$$
where $\mu_0$ is given by Lemma \ref{lmq5.2}. Then $\mu_*:=\sup\Sigma_*\leq\mu^*$ and $(0,\mu_*)\subset\Sigma_*$. Similar to the above, it can be obtained that $\mu_*\in\Sigma_*$. The proof is completed.
\ \ \ \ \fbox{}

\section{Discussion}
\setcounter{equation}{0}

In this paper, we have examined a predator-prey model with double free boundaries
$x= g(t)$ and $x=h(t)$ for the predator, which describes the movement process through
the two free boundaries. The dynamic behavior are discussed.

A great deal of previous mathematical investigation on the spreading of population has been based on the traveling wave fronts of the predator-prey system over the entire space $\mathbb{R}$
 \bes
 \left\{\begin{array}{lll}
 u_t-u_{xx}=u(1-u+av), &t>0,\ \ x\in\mathbb{R},\\[1mm]
  v_t-Dv_{xx}=v(b-v-cu),\ \ &t>0, \ \ x\in\mathbb{R}.
 \end{array}\right.\label{6.1}
 \ees
A striking difference between (\ref{1.2}) and (\ref{6.1}) is that the spreading front in (\ref{1.2}) is given explicitly by a function $x=h(t)$, beyond which the population density of the predator is $0$, while in (\ref{6.1}), the population $u(t, x)$ becomes positive for all $x$ once $t$ is positive. Second, (\ref{6.1}) guarantees successful spreading of the predator species for any nontrivial initial population $u(0, x)$, regardless of its initial size and
supporting area, but the dynamics of (\ref{1.2}) exhibits a spreading-vanishing dichotomy. The phenomenon exhibited by this dichotomy seems closer to the reality.

The spreading-vanishing dichotomy results also indicate that:

(i) When the spreading happens, both the predator and prey will converge to positive
constants for the weakly hunting case, while
the predator will converge to a positive constant and the prey will vanish for the strongly hunting case. These dynamic behaviours are similar to that of solution to the Cauchy problem
of (\ref{6.1}).

(ii) When the vanishing occurs, the predator will vanish and the prey will converge to a positive constant.

The criteria governing spreading and vanishing (Theorems \ref{th5.1} and \ref{th5.2}) tell us that whether spreading or vanishing are completely determined by sizes of both the initial habitat and initial data of the predator, and the moving parameter/coefficient $\mu$ of free boundaries.

These results tell us that in order to control the prey species (pest species) we should
put predator species (natural enemies) at the initial state at least in one of
three ways: (i) expand the predator's targets, (ii) increase the
moving parameter/coefficient of free boundaries, (iii) augment the initial
density of the predator species.

\vskip 25pt\begin{appendix}
\noindent{\bf\LARGE Appendix}


\section{Global estimate of the solution $w$ to (\ref{4.3})}
 \def\theequation{\Alph{section}.\arabic{equation}}
 \setcounter{equation}{0}

\begin{prop}\label{pa.1} Let $(u,v,g,h)$ be any solution of $(\ref{1.2})$ and
assume $h_\infty-g_\infty<\infty$, if $w(t,y)$ is the solution of $(\ref{4.3})$,
then there exists a constant $K_0$ such that
\begin{equation}
\|w\|_{C^{\frac{1+\alpha}{2},1+\alpha}([1,\infty)\times[-1,1])}<K_0.\label{A.1}
\end{equation}\end{prop}

{\bf Proof.} We are inspired by \cite[Theorem A2]{BKD}. For convenience,
we denote $\varphi_n(t)=\varphi(t+n), \psi_n=\psi(t+n,y), z_n=z(t+n,y), w_n=w(t+n,y)$.
Let $w(t+n,y)=a^n(t,y)+b^n(t,y)$, where $a^n$ and $b^n$ are solutions of
 \bess
\left\{\begin{array}{ll}
a^n_t=\varphi_n(t)a^n_{yy}, \ \ &t>0,-1<y<1,\\[1mm]
a^n(t,-1)=a^n(t,1)=0,\ \ &t>0,\\[1mm]
a^n(0,y)=w(n,y),&-1\leq y\leq 1
\end{array}\right.\eess
 and
  \bess
\left\{\begin{array}{ll}
b^n_t=\varphi_n(t)b^n_{yy}+\psi_nb^n_y+\psi_na^n_y+w_n(1-w_n+az_n), \ \
&t>0,-1<y<1,\\[1mm]
b^n(t,-1)=b^n(t,1)=0, &t>0,\\[1mm]
b^n(0,y)=0,&-1\leq y\leq 1,
\end{array}\right.
  \eess
respectively.  Let $0<\lambda_1\leq\lambda_2\leq\cdots$ be the eigenvalues of the
problem
 $$-\phi_{yy}=\lambda\phi,\ -1<y<1;\ \ \phi(-1)=\phi(1)=0,$$
and let $\phi_1,\phi_2,\cdots$ be the corresponding set of orthonormal eigenfunctions.
We may express $a^n$ as
 \[a^n(t,y)=\sum\limits_{k\geq
1}\exp\kk(-\lambda_k\int_0^t\varphi_n(s)ds\rr)w^n_k\phi_k,\]
where $w^n_k=(\phi_k,w(n,\cdot))_{L^2}$. In view of $0<u\leq M_1$ (cf. Lemma
\ref{lm2.1}), it follows that
 \[\begin{array}{ll}
\kk\|\dd\frac{\partial^{2j}a^n(T)}{\partial
y^{2j}}\rr\|_{L^2(-1,1)}^2&=\sum\limits_{k\geq
1}\lambda_k^{2j}\exp\dd\kk(-2\lambda_k\int_0^T\varphi_n(s)ds\rr)(w^n_k)^2\\[3mm]
&\leq\sup\limits_{\ell\geq 0}\left\{\ell^{2j}\exp\kk(\frac{-8T\ell}
{(h_\infty-g_\infty)^2}\rr)\rr\}\|w(n,x)\|_{L^2(-1,1)}^2\\[3mm]
&\leq2M_1^2\kk(\dd\frac{(h_\infty-g_\infty)^2j}{4T}\rr)^{2j}{\rm e}^{-2j}.
 \end{array}\]
By the $L^p$ estimates and Sobolev's imbedding theorem, we therefore have that
$\|a^n(t)\|_{C^2[-1,1]}\leq K_1(1+t^{-j})$, where $K_1$ is independent of $n$
provided that $j\geq 2$. From this last estimate and the differential equation
satisfied by $a^n$, we obtain
$\|a^n_t(t)\|_{C[-1,1]}\leq K_1h_0^{-2}(1+t^{-j})$.
Hence, there exists positive constant $K_2$ such that
$\|a^n\|_{C^{\frac{1+\alpha}{2},1+\alpha}(E_1)}<K_2$, where $E_1=[{1\over
2},2]\times[-1,1]$ and $K_2$ depends only on $K_1$ and $E_1$.

Next we estimate $b^n$. It is obvious that the function $c^n={\rm e}^{-\frac{1}{t}}b^n$ satisfies
  $$\left\{\begin{array}{ll}
c^n_t=\varphi_n(t)c^n_{yy}+\psi_nc^n_y+f_n(t,x), \ \ &t>0,-1<y<1,\\[1mm]
b^n(t,-1)=b^n(t,1)=0, &t>0,\\[1mm]
b^n(0,y)=0,&-1\leq y\leq 1,
 \end{array}\right.$$
where
  $$f_n(t,x)=\left\{\begin{array}{ll}
\dd\frac{1}{t^2}{\rm e}^{-\frac{1}{t}}w_n+(\psi_n+\dd\frac{1}{t^2})
{\rm e}^{-\frac{1}{t}}a^n_y+{\rm e}^{-\frac{1}{t}}w_n(1-w_n+az_n),&\ t>0\\[2mm]
0,&\ t=0.
  \end{array}\right.$$
Note that $\lim\limits_{t\to 0^+}\dd t^{-j}{\rm e}^{-\frac{1}{t}}=0$ for any $j>0$, we have
$f_n(t,x)$ is continuous in $E=[0,3]\times[-1,1]$
and $\|f_n\|_{C(E)}\leq K_3$ where $K_3$ is dependent on $K_1$ and independent of $n$.
By using of \cite[Theorem 4, p191]{Fri}, we can obtain that
$\|c^n\|_{C^{\frac{1+\alpha}{2},1+\alpha}(E_1)}<\tilde K_3$ where
$\tilde K_3$ depends on $K_3$. It follows that
$\|b^n\|_{C^{\frac{1+\alpha}{2},1+\alpha}(E_1)}<K_4$. We therefore have that
 \[\|w\|_{C^{\frac{1+\alpha}{2},1+\alpha}(E_n)}\leq
\|a^n\|_{C^{\frac{1+\alpha}{2},1+\alpha}(E_1)}+\|b^n\|_{C^{\frac{1+\alpha}{2},1+\alpha}(E_1)}
 <K_2+K_4=K_0,\]
where $E_n=[n+\frac{1}{2},n+2]\times[-1,1]$. It easily to get (\ref{A.1}) since the
intervals $E_n$ overlap and $K_0$ is independent of $n$.\ \ \ \ \fbox{}

\section{Estimates of solutions to parabolic partial differential inequalities}
   \def\theequation{\Alph{section}.\arabic{equation}}
 \setcounter{equation}{0}

Let $d,\beta$ and $\theta$ be fixed positive constants. In order to investigate the long time behavior of the solution $(u,v)$ to (\ref{1.2}), we should prove the following two propositions.

\begin{prop}\label{pB.1} \ For any given $\ep>0$ and  $L>0$, there exist
$l_\ep>\max\big\{L,\frac{\pi}2\sqrt{d/\beta}\big\}$ and $T_\ep>0$, such that when the
continuous and non-negative function $w(t,x)$ satisfies
 $$ \left\{\begin{array}{ll}
  w_t-dw_{xx}\geq\, (\leq)\, w(\beta-\theta w), \ \ &t>0, \ \ -l_\ep<x<l_\ep,\\[2.5mm]
 w(t,\pm l_\ep)\geq\, (=)\, 0,&t\ge 0,
 \end{array}\right.$$
and $w(0,x)>0$ in $(-l_\ep,l_\ep)$, then
\bess w(t,x)>\beta/\theta-\ep \ \ \big(w(t,x)<\beta/\theta+\ep\big), \ \ \
  \forall \ t\geq T_\ep, \ \ x\in[-L,L].\eess
Which implies
 \[\liminf_{t\to\infty}w(t,x)>\beta/\theta-\ep \ \ \kk(\limsup_{t\to\infty}w(t,x)<\beta/\theta+\ep\rr)
 \ \ \ \mbox{uniformly\, on }\, [-L,L].\]
  \end{prop}

{\bf Proof.} \ Let $l>\frac{\pi}2\sqrt{d/\beta}$ be a parameter. Assume that $w_l(x)$
is the unique positive solution of
  \bes \left\{\begin{array}{ll}
  -dw_{xx}=w(\beta-\theta w),\ \ -l<x<l,\\[1mm]
 w(\pm l)=0.
 \end{array}\right.\label{B.1}\ees
By Lemma 2.2 of \cite{DM}, $\lim_{l\to\infty}w_l(x)=\beta/\theta$ uniformly in any compact subset of $\mathbb{R}$. So, for any given $L>0$ and $\ep>0$, there exists $l_\ep>\max\big\{L,\frac{\pi}2\sqrt{d/\beta}\big\}$, which also
depends on $d,\beta$ and $\theta$, such that
 \bes \beta/\theta-\ep/2<w_l(x)<\beta/\theta+\ep/2, \ \ \ \forall \ l\geq l_\ep, \
x\in[-L,L].\label{B.2}\ees

Let $w_0(x)\in C([-l_\ep,l_\ep])$ be a positive function and $w_\ep(t,x)$ be the
unique solution of
  $$ \left\{\begin{array}{ll}
  w_{t}-dw_{xx}=w(\beta-\theta wp), \ \ &t>0, \ \ -l_\ep<x<l_\ep,\\[1mm]
 w(t,\pm l_\ep)=0,&t\ge 0,\\[1mm]
 w(0,x)=w_0(x),&-l_\ep\leq x\leq l_\ep.
 \end{array}\right.$$
Since $l_\ep>\frac{\pi}2\sqrt{d/\beta}$, it is well known that $\lim_{t\to\infty}
w_\ep(t,x)=w_{l_\ep}(x)$ uniformly in the compact subset of $(-l_\ep,l_\ep)$. Thanks
to (\ref{B.2}), there is a $T_\ep\gg 1$ such that
 \bes \beta/\theta-\ep<w_\ep(t,x)<\beta/\theta+\ep, \ \ \
  \forall \ t\geq T_\ep, \ \ x\in[-L,L].\label{B.3}\ees
Our conclusion is followed from (\ref{B.3}) and the comparison principle.\ \ \ \ \fbox{}

\begin{prop}\label{pB.2} \ Let $k$ be a positive constant. For any given $\ep>0$ and $L>0$,
there exist $l_\ep>\max\big\{L,\frac{\pi}2\sqrt{d/\beta}\big\}$ and $T_\ep>0$,
such that when the continuous and non-negative function $z(t,x)$ satisfies
 $$ \left\{\begin{array}{ll}
  z_t-dz_{xx}\geq\, (\leq)\, z(\beta-\theta z), \ \ &t>0, \ \ -l_\ep<x<l_\ep,\\[1mm]
 z(t,\pm l_\ep)\geq\, (\leq)\, k,&t\ge 0,
 \end{array}\right.$$
and $z(0,x)>0$ in $(-l_\ep,l_\ep)$, then we have
 \[ z(t,x)>\beta/\theta-\ep \ \ \big(z(t,x)<\beta/\theta+\ep\big), \ \ \ \forall \ t\geq T_\ep, \ \ x\in[-L,L].\]
This implies
 \[\liminf_{t\to\infty}z(t,x)\geq \beta/\theta-\ep \ \ \kk(\limsup_{t\to\infty}z(t,x)<\beta/\theta+\ep\rr)
  \ \ \ \mbox{uniformly\, on }\, [-L,L].\]
 \end{prop}

{\bf Proof.} Let $l>\frac{\pi}2\sqrt{d/\beta}$ be a parameter, and $z_l(x)$ the
unique positive solution of
 \bes \left\{\begin{array}{ll}
  -dz_{xx}=z(\beta-\theta z),\ \ -l<x<l,\\[1mm]
 z(\pm l)=k.
 \end{array}\right.\label{B.4}\ees
(refer to the proof of Lemma 2.3 in \cite{DM}). We claim that
 \bes
  \lim_{l\to\infty}z_l(x)=\beta/\theta \ \ \mbox{ uniformly\ in\ any\ compact\ subset\ of} \, \  \mathbb{R}.\label{B.5}\ees

For the case $k>\beta/\theta$. By the maximum principle
we see that $\beta/\theta\leq z_l(x)\leq k$ for all $x\in[-l,l]$. Note that $z_l(x)\leq k$, by
the comparison principle we have that $z_l(x)$ is decreasing in $l$. Therefore, the
limit $\lim_{l\to \infty}z_l(x)=z(x)$ exists, and $z(x)\geq \beta/\theta$ and $z(x)$
satisfies
 \[-dz_{xx}=z(\beta-\theta z),\ \ x\in\mathbb{R}.\]
By Theorem 1.2 of \cite{DM}, $z(x)\equiv \beta/\theta$. Using the interior estimate we
assert that $\lim_{l\to \infty}z_l(x)=z(x)$ uniformly in any compact subset of
$\mathbb{R}$. Hence (\ref{B.5}) holds.

For the case $k\leq \beta/\theta$. Choose $k_0>\beta/\theta$ and let $z_l^0$ be the
unique positive solution of (\ref{B.4}) with $k=k_0$. By the comparison principle we
have $w_l(x)\leq z_l(x)\leq z_l^0(x)$ in $[-l,l]$, where $w_l(x)$ is the unique
positive solution of (\ref{B.1}) with $l>\frac{\pi}2\sqrt{d/\beta}$. Take into account the result we have proved in the above and Lemma 2.2 of \cite{DM}, it is deduced that
(\ref{B.5}) holds.

In view of (\ref{B.5}), for any given $L>0$ and  $\ep>0$, there is
$l_\ep>\max\big\{L,\frac{\pi}2\sqrt{d/\beta}\big\}$, which also depends on $d,\beta,\theta$ and $k$,  such that
 \bes \beta/\theta-\ep/2<z_l(x)<\beta/\theta+\ep/2, \ \ \ \forall\ l\geq l_\ep, \
 x\in[-L,L].\label{B.6}\ees
Let $z_0(x)\in C([-l_\ep,l_\ep])$ be a positive function and $z_\ep(t,x)$ be the
solution of
  $$ \left\{\begin{array}{ll}
  z_{t}-dz_{xx}=z(\beta-\theta z), \ \ &t>0, \ \ -l_\ep<x<l_\ep,\\[1mm]
 z(t,\pm l_\ep)=k,&t\ge 0,\\[1mm]
 z(0,x)=z_0(x),&-l_\ep\leq x\leq l_\ep.
 \end{array}\right.$$
Recall $l_\ep>\frac{\pi}2\sqrt{d/\beta}$, we shall illustrate that
 \bes\dd\lim_{t\to\infty} z_\ep(t,x)=z_{l_\ep}(x)\ \ \mbox{ uniformly\ in\ the\
compact\ subset\ of} \ \ (-l_\ep,l_\ep).
 \label{B.7}\ees
In fact, take a positive constant $q$ and let $\phi_q(t,x)$ be the unique solution of
  \bess \left\{\begin{array}{ll}
  \phi_t-d \phi_{xx}=\phi(\beta-\theta \phi), \ \ &t>0, \ \ -l_\ep<x<l_\ep,\\[1mm]
 \phi(t,\pm l_\ep)=k,&t\ge 0,\\[1mm]
 \phi(0,x)=q,&-l_\ep\leq x\leq l_\ep.
 \end{array}\right.\eess
Let $M\gg 1$ and $0<m\ll 1$. Then $M$ and $m$ are the ordered upper and lower solutions of (\ref{B.4}) with $l=l_\ep$. Therefore,
$\phi_M(t,x)$ is monotone decreasing and $\phi_m(t,x)$ is monotone increasing in $t$. So, the limits
$\lim_{t\to\infty}\phi_M(t,x)=\phi_M(x)$ and $\lim_{t\to\infty}\phi_m(t,x)=\phi_m(x)$ exist,  and they are all positive solution of (\ref{B.4}) with $l=l_\ep$. Hence, $\phi_M(x)=\phi_m(x)=z_{l_\ep}(x)$. Meanwhile, the comparison
principle yields $\phi_m(t,x)\leq z_\ep(t,x)\leq \phi_M(t,x)$. Consequently,
$\lim_{t\to\infty} z_\ep(t,x)=z_{l_\ep}(x)$. By use of the interior estimate, it can be shown that this limit is uniformly in the compact subset of $(-l_\ep,l_\ep)$.

Thanks to (\ref{B.6}) and (\ref{B.7}), there is $T_\ep\gg 1$ such that
 \bess \beta/\theta-\ep<z_\ep(t,x)<\beta/\theta+\ep, \ \ \
  \forall \ t\geq T_\ep, \ \ x\in[-L,L].\eess
By sue of this fact and the comparison principle, the proof is immediately completed. \ \ \ \ \fbox{}
\end{appendix}

\vskip 10pt
{\bf Acknowledgement}\ The authors would like to thank Professor
Xinfu Chen for several useful discussions about this subject.

\end{document}